\documentclass[a4paper,12pt,twoside]{article}
\usepackage{amsmath,amssymb,ifthen}
\usepackage[amsmath,thmmarks]{ntheorem}
\usepackage{paper}
\usepackage[all]{xy}
\usepackage{mathrsfs}
\SelectTips{cm}{11}
\newcommand{\rig}{\mathrm{rig}}
\newcommand{\ad}{\mathrm{ad}}

\DeclareMathOperator{\Ind}{Ind}
\DeclareMathOperator{\Spf}{Spf}
\DeclareMathOperator{\Spa}{Spa}
\DeclareMathOperator{\rank}{rank}

\DeclareMathOperator{\spp}{sp}

\DeclareMathOperator{\Sh}{Sh}

\DeclareMathOperator{\Res}{Res}
\DeclareMathOperator{\Lie}{Lie}
\DeclareMathOperator{\Spl}{Spl}

\title{Compactly supported cohomology and nearby cycle cohomology of open Shimura varieties of PEL type}
\author{Naoki Imai and Yoichi Mieda}
\def\headertitle{Cohomology of open Shimura varieties of PEL type}

\begin{document}

\maketitle

\begin{firstfootnote}
\end{firstfootnote}

\begin{abstract}
 In this paper, we compare two cohomology groups associated to Shimura varieties
 of PEL type, which are not proper over the base. One is the compactly supported
 $\ell$-adic cohomology, and the other is the nearby cycle cohomology,
 namely, the compactly supported cohomology of
 the nearby cycle complex for the canonical integral model of the Shimura variety
 over $\Z_p$.  
 We prove that the $G(\Q_p)$-cuspidal part of these cohomology groups are
 the same, where $G$ denotes the reductive algebraic group naturally attached to the PEL datum.
 Some applications to unitary Shimura varieties are also given.
\end{abstract}

\section{Introduction}

Study of the $\ell$-adic cohomology of Shimura varieties is an important theme in number theory and
arithmetic algebraic geometry. Let $\mathrm{Sh}_K=\Sh_K(G,X)$ be a Shimura variety associated to
a Shimura datum $(G,X)$ and a compact open subgroup $K$ of $G(\A^\infty)$.
As is well-known, it is naturally defined over a number field $E$, called the reflex field,
which is canonically attached to $(G,X)$. 

First of all, let us consider the case where our Shimura variety is compact; namely, $\Sh_K$ is proper over $E$.
For an irreducible algebraic representation $\xi$ of $G$ over $\overline{\Q}_\ell$, 
the $\ell$-adic \'etale cohomology
\[
 H^i(\Sh,\mathcal{L}_\xi)=\varinjlim_{K}H^i(\Sh_K\otimes_E\overline{E},\mathcal{L}_\xi)
\]
of the associated local system $\mathcal{L}_\xi$ is naturally equipped with the action of
$G(\A^\infty)\times\Gal(\overline{E}/E)$. 
By Matsushima's formula, this cohomology group is strongly related to the automorphic spectrum of $G(\A)$.
If moreover our Shimura variety is so-called PEL type, namely, obtained as a moduli space of abelian varieties
with several additional structures, then it has a natural integral model over $\mathcal{O}_{E_v}$, the ring of
integers of the completion of $E$ at a finite place $v$. This integral model, also denoted by $\Sh_K$,
is very useful to investigate the action of the Galois group $\Gal(\overline{E}/E)$ on $H^i(\Sh,\mathcal{L}_\xi)$.
If $\Sh_K$ is smooth over $\mathcal{O}_{E_v}$, the proper smooth base change theorem says that
the $\ell$-adic cohomology $H^i(\Sh_K\otimes_E\overline{E},\mathcal{L}_\xi)$ is canonically isomorphic to
the $\ell$-adic cohomology of the special fiber $\Sh_{K,v}$ of the integral model. As a result,
the $\Gal(\overline{E}/E)$-representation $H^i(\Sh_K\otimes_E\overline{E},\mathcal{L}_\xi)$ is unramified at
$v$, and the alternating sum of the traces of $\Frob_v$, the Frobenius element at $v$,
on $H^i(\Sh_K\otimes_E\overline{E},\mathcal{L}_\xi)$ can be computed by the Lefschetz trace formula.
Even if $\Sh_K$ is not smooth over $\mathcal{O}_{E_v}$, we can investigate the action of $\Gal(\overline{E}_v/E_v)$
on $H^i(\Sh_K\otimes_E\overline{E},\mathcal{L}_\xi)$ by using the theory of nearby cycle functor.
Such a method provides an efficient way to establish a part of the Langlands correspondence;
see \cite{MR1876802}, \cite{MR2074714}, \cite{MR2074715}, \cite{MR2169874}, \cite{Shin-Annals}
and \cite{Scholze-LLC}, for instance.
The crucial point is that the properness of $\Sh_K$ over $\mathcal{O}_{E_v}$ ensures the coincidence of
the cohomology of the generic fiber and the nearby cycle cohomology, that is, 
the cohomology of the nearby cycle complex over the special fiber. 

Now we pass to the case where $\Sh_K$ is not compact, while is of PEL type.
There are many interesting cases satisfying these conditions, including the modular curves
and the Siegel modular varieties. For the $\ell$-adic cohomology of a Shimura variety over $E$,
we have at least three choices; the ordinary $\ell$-adic cohomology,
the compactly supported $\ell$-adic cohomology and the intersection cohomology.
Usually we prefer the intersection cohomology, since it has a close relation to the automorphic
spectrum of $G(\A)$. However, the theory of minimal compactifications tells us that
these three cohomology groups are the same ``up to induced representations''.

Can we use the natural integral model in this case to investigate the Galois action on the $\ell$-adic cohomology?
If $\Sh_K$ is smooth over $\mathcal{O}_{E_v}$, Fujiwara's trace formula \cite{MR1431137} is a very powerful tool
to carry it out. By using it, one can compute the alternating sum of the traces of the action of $\Frob_v$
on the intersection cohomology $\mathit{IH}^i(\Sh_K\otimes_{E}\overline{E},\mathcal{L}_\xi)$;
see \cite{MR2350050}, \cite{MR2567740}. On the other hand, if $\Sh_K$ has bad reduction,
we cannot use exactly the same method as in the compact case, since the nearby cycle cohomology
is different from the $\ell$-adic cohomology of the generic fiber, whichever cohomology we choose among
the three mentioned above.

The purpose of this paper is to measure the gap between the compactly supported cohomology (of the Shimura
variety over $E$) and the nearby cycle cohomology, that is, the compactly supported cohomology of
the nearby cycle complex over the special fiber of the integral model.
Our main theorem claims that the gap is non-cuspidal at $p$, where $p$ is the prime divisible by $v$.
More precise statement is the following:

\begin{thm}[Theorem \ref{thm:main-thm}]\label{thm:main-thm-intro}
 In the kernel and the cokernel of the natural map
 \[
 \varinjlim_{K} H^i_c(\Sh_{K,\overline{v}},R\psi\mathcal{L}_\xi)\longrightarrow 
 \varinjlim_{K} H^i_c(\Sh_K\otimes_{E_v}\overline{E}_v,\mathcal{L}_\xi),
 \]
 no irreducible supercuspidal representation of $G(\Q_p)$ appears as a subquotient.
\end{thm}
Recall that an irreducible smooth representation of $G(\Q_p)$ is said to be supercuspidal
if it does not appear as a subquotient of parabolically induced representations from any proper parabolic subgroup.

In particular, for an irreducible admissible representation $\Pi$ of $G(\A^\infty)$ such that
$\Pi_p$ is supercuspidal, we have an isomorphism of $\Gal(\overline{E}_v/E_v)$-modules
\begin{equation*}
 \Bigl(\varinjlim_{K} H^i_c(\Sh_{K,\overline{v}},R\psi\mathcal{L}_\xi)\Bigr)[\Pi]\cong 
 \Bigl(\varinjlim_{K} H^i_c(\Sh_K\otimes_{E_v}\overline{E}_v,\mathcal{L}_\xi)\Bigr)[\Pi].
 \tag{$*$}
\end{equation*}
Here $(-)[\Pi]$ denotes the $\Pi$-isotypic part (for precise definition, see Notation below).
As a result, for a representation $\Pi$ as above, the contribution of it to the compactly supported cohomology
of the Shimura variety can be investigated by using the nearby cycle functor.
Note that the same thing holds for the ordinary cohomology and the intersection cohomology,
since the gaps between the compactly supported cohomology and these two cohomology groups are non-cuspidal,
as already explained.

If an irreducible admissible representation $\Pi$ of $G(\A^\infty)$ has a supercuspidal local component
$\Pi_{p'}$ at a prime $p'\neq p$, then it is easy to obtain an isomorphism
 $(*)$ by using the minimal compactification
with hyperspecial level at $p$.
The main difficulty to prove Theorem \ref{thm:main-thm-intro} is the absence of a good compactification
of our integral Shimura variety with higher level at $p$. 
Our main idea to overcome it is to use rigid geometry. Roughly speaking, the gap between
the compactly supported cohomology and the nearby cycle cohomology can be measured by the rigid subspace
of the generic fiber of our integral model consisting of points corresponding to
abelian varieties with bad reduction. We will find a nice stratification on this rigid subspace
to conclude that it is ``geometrically'' induced from proper parabolic subgroups of $G(\Q_p)$.
To construct the stratification, we use the theory of arithmetic toroidal compactifications
with hyperspecial level at $p$. Such compactifications were constructed by Faltings and Chai \cite{MR1083353}
for Siegel modular varieties, and by Lan \cite{Kai-Wen} for general Shimura varieties of PEL type.
There is also a work by Fujiwara on this topic (\cf \cite[\S 7]{MR2310255}).

It will be worth noting that our method is totally geometric and thus is also valid
in the torsion coefficient case.
See Theorem \ref{thm:main-thm-mod-l} for an analogue of Theorem \ref{thm:main-thm-intro} with
the $\overline{\F}_\ell$-coefficients. 

We sketch the outline of this paper. In Section 2, we give some preliminary results for constructing
a nice stratification.
In Section 3, after recalling some notation concerning with Shimura varieties of PEL type, we give
a key geometric construction; namely, the construction of a stratification on the rigid subspace
that measures the gap of two cohomology groups. 
We also observe that this stratification behaves well under the Hecke action on the tower of Shimura varieties.
In Section 4, we introduce several cohomology groups of the Shimura variety and prove 
Theorem \ref{thm:main-thm-intro}.
In Section 5, we give some applications of our main result to unitary Shimura varieties.
The authors expect that our result has more interesting applications in a symplectic case.
We will study it in our future work.

\bigbreak

\noindent{\bfseries Notation}\quad
Put $\widehat{\Z}=\prod_{\text{$p$: prime}}\Z_p$ and $\A^\infty=\widehat{\Z}\otimes_\Z\Q$.
For a prime $p$, put $\widehat{\Z}^p=\prod_{\text{$p'\neq p$: prime}}\Z_{p'}$ and 
$\A^{\infty,p}=\widehat{\Z}^p\otimes_\Z\Q$.

Assume that we are given two locally profinite groups $G$ and $\Gamma$ and a prime $\ell$.
For an admissible/continuous representation $V$ of $G\times \Gamma$ over $\overline{\Q}_\ell$ in the sense of \cite[I.2]{MR1876802} and
an irreducible admissible representation $\pi$ of $G$, put $V[\pi]=\bigoplus_{\sigma}\sigma^{\oplus m_{\pi\boxtimes\sigma}}$,
where $\sigma$ runs through finite-dimensional irreducible continuous $\overline{\Q}_\ell$-representations of $\Gamma$,
and $m_{\pi\boxtimes\sigma}$ denotes the coefficient of $[\pi\boxtimes\sigma]$
in the image of $V$ in the Grothendieck group $\mathrm{Groth}_{G\times\Gamma,\ell}(G\times\Gamma)$ considered in \cite[I.2]{MR1876802}.
It is a semisimple continuous representation of $\Gamma$.

Every sheaf and cohomology are considered in the \'etale topology.

\section{Rigid geometry and semi-abelian schemes}\label{sec:rig-semi-ab}
\subsection{Notation for adic spaces}\label{subsec:notation-adic}
Throughout this paper, we will use the framework of adic spaces introduced by Huber
(\cf \cite{MR1207303}, \cite{MR1306024}, \cite{MR1734903}). Here we recall some notation briefly.

Let $S$ be a noetherian scheme and $S_0$ a closed subscheme of $S$. We denote the formal completion
of $S$ along $S_0$ by $\mathcal{S}$. Put $\mathcal{S}^\rig=t(\mathcal{S})_a$, where
$t(\mathcal{S})$ is the adic space associated to $\mathcal{S}$ (\cf \cite[\S 4]{MR1306024})
and $t(\mathcal{S})_a$ denotes the open adic subspace of $t(\mathcal{S})$
consisting of analytic points.
It is a quasi-compact analytic adic space.

Let $X$ be a scheme of finite type over $S$. Put $X_0=X\times_SS_0$ and denote
the formal completion of $X$ along $X_0$ by $\widehat{X}$. Then we can construct an adic space
$\widehat{X}^\rig$ in the same way as $\mathcal{S}^\rig$. The induced morphism
$\widehat{X}^\rig\longrightarrow \mathcal{S}^\rig$ is of finite type.
On the other hand, we can construct another adic space $X\times_S\mathcal{S}^\rig$.
 Indeed, since we have morphisms of locally ringed spaces
$(\mathcal{S}^\rig,\mathcal{O}_{\mathcal{S}^\rig})\longrightarrow (t(\mathcal{S}),\mathcal{O}_{t(\mathcal{S})})\longrightarrow (S,\mathcal{O}_S)$ (for the second one, see \cite[Remark 4.6 iv)]{MR1306024}), we can make the fiber product $X\times_S\mathcal{S}^\rig$ in the sense of
\cite[Proposition 3.8]{MR1306024}. 
For simplicity, we write $X^\ad$ for $X\times_S\mathcal{S}^\rig$, though it depends on
$(S,S_0)$.
Since the morphism $\mathcal{S}^\rig\longrightarrow S$ factors through $S^0=S\setminus S_0$,
we have $X\times_S\mathcal{S}^\rig=(X\times_SS^0)\times_{S^0}\mathcal{S}^\rig$.
In particular, $X^\ad$ depends only on $X\times_SS^0$.
The natural morphism $X^\ad\longrightarrow \mathcal{S}^\rig$ is locally of finite type,
but not necessarily quasi-compact; see the following example.

\begin{exa}
 Let $R$ be a complete discrete valuation ring and $F$ its fraction field.
 Consider the case where $S=\Spec R$ and $S_0$ is defined by the maximal ideal of $R$.
 Then, for an $S$-scheme $X$ of finite type, $X^\ad$ can be regarded as the rigid space over $F$
 associated to a scheme $X\times_S\Spec F$ over $F$.
 For example, $(\A^1_S)^\ad=(\A^1_F)^\ad$ is the rigid-analytic affine line over $F$ and thus is not
 quasi-compact. On the other hand, $(\widehat{\A}^1_S)^\rig$ is the unit disc
 ``$\lvert z\rvert\le 1$'' in $(\A^1_F)^\ad$, which is quasi-compact.
\end{exa}

\begin{lem}\label{lem:fiber-product}
 The functors $X\longmapsto \widehat{X}^\rig$ and $X\longmapsto X^\ad$ commute with
 fiber products.
\end{lem}

\begin{prf}
 For the functor $X\longmapsto \widehat{X}^\rig$, it can be checked easily
 (\cf \cite[Lemma 3.4]{MR2211156} and \cite[Proof of Lemma 4.3 v)]{formalnearby}).
 Consider the functor $X\longmapsto X^\ad$. 
 Let $Y\longrightarrow X\longleftarrow Z$ be a diagram of $S$-schemes of finite type.
 What we should prove is
 \[
  (Y\times_XZ)\times_S\mathcal{S}^\rig\cong (Y\times_S\mathcal{S}^\rig)\times_{X\times_S\mathcal{S}^\rig}(Z\times_S\mathcal{S}^\rig).
 \]
 It is not totally automatic, since $Y\times_XZ$ is not a fiber product in the category of locally ringed spaces.
 It follows from the fact that morphisms of locally ringed spaces
 $\Spa (A,A^+)\longrightarrow \Spec B$ correspond bijectively to ring homomorphisms $B\longrightarrow A$
 (this fact is used implicitly in \cite[Remark 4.6 (iv)]{MR1306024} to define $t(\mathcal{S})\longrightarrow S$).
\end{prf}

Let us compare $\widehat{X}^\rig$ and $X^\ad$; by the commutative diagram
\[
 \xymatrix{%
 \widehat{X}^\rig\ar[r]\ar[d]& X\ar[d]\\
 \mathcal{S}^\rig\ar[r]& S
 }
\]
and the universality of the fiber product $X\times_S\mathcal{S}^\rig$, we have a natural
morphism $\widehat{X}^\rig\longrightarrow X^\ad$.

\begin{lem}\label{lem:rig-ad}
 \begin{enumerate}
  \item If $X$ is separated over $S$, the natural morphism $\widehat{X}^\rig\longrightarrow X^\ad$
	is an open immersion.
  \item If $X$ is proper over $S$, the natural morphism $\widehat{X}^\rig\longrightarrow X^\ad$
	is an isomorphism.
 \end{enumerate}
\end{lem}

\begin{prf}
 See \cite[Remark 4.6 (iv)]{MR1306024}.
\end{prf}

\begin{rem}\label{rem:rig-ad-base-change}
 Let $f\colon S'\longrightarrow S$ be a morphism of finite type and $S'_0=S'\times_SS_0$.
 We denote by $\mathcal{S}'^\rig$ the formal completion of $S'$ along $S'_0$.
 Then, all constructions above are compatible with the base change
 by $f$. More precisely, for a scheme $X$ of finite type over $S$, 
 we have $(X\times_SS')^{\wedge\rig}\cong \widehat{X}^\rig\times_{\mathcal{S}^\rig}\mathcal{S}'^\rig$
 and $(X\times_SS')^{\text{$S'$-ad}}\cong X^\ad\times_{\mathcal{S}^\rig}\mathcal{S}'^\rig$.
 Here $(-)^{\text{$S'$-ad}}$ denotes the functor $(-)^\ad$ for the base $(S',S'_0)$,
 namely, $(-)^\text{$S'$-ad}=(-)\times_{S'}\mathcal{S}'^\rig$.
\end{rem}

The following notation is used in Section \ref{sec:Sh-var}. In the remaining part of this subsection,
assume that $S$ is the spectrum of a complete discrete valuation ring and
$S_0$ is the closed point of $S$. For a scheme of finite type $X$ over $S$, 
we have a natural morphism of locally and topologically ringed spaces 
$(t(\widehat{X}),\mathcal{O}_{\!t(\widehat{X})}^+)\longrightarrow (\widehat{X},\mathcal{O}_{\widehat{X}})$
(\cf \cite[Proposition 4.1]{MR1306024}).
Note that the underlying continuous map $t(\widehat{X})\longrightarrow X_0$ is different from the map
$t(\widehat{X})\longrightarrow X$ considered above.
We denote the composite $\widehat{X}^\rig\hooklongrightarrow t(\widehat{X})\longrightarrow X_0$
by $\spp_{\widehat{X}}$, or simply by $\spp$.

Let $Y$ be a closed subscheme of $X_0$ and $\mathcal{X}$ the formal completion of $X$ along $Y$.
Then we can consider the generic fiber $t(\mathcal{X})_\eta=S^0\times_St(\mathcal{X})$ of the adic space $t(\mathcal{X})$.
This is so-called the rigid generic fiber of $\mathcal{X}$ due to Raynaud and Berthelot, in the context of
adic spaces.
If $Y=X_0$, then $t(\mathcal{X})_\eta=\widehat{X}^\rig$.

\begin{lem}\label{lem:rig-gen-fiber}
 The natural morphism $t(\mathcal{X})_\eta\longrightarrow \widehat{X}^\rig$ induced from
 $\mathcal{X}\longrightarrow \widehat{X}$ is an open immersion.
 Its image coincides with $\spp^{-1}(Y)^\circ$, where $(-)^\circ$ denotes the interior (in $\widehat{X}^\rig$).
\end{lem}

\begin{prf}
 See \cite[Lemma 3.13 i)]{MR1620118}.
\end{prf}

\subsection{Etale sheaves associated to semi-abelian schemes}\label{subsec:etale-sheaves-semi-abelian}
We continue to use the notation introduced in the beginning of the previous subsection.
Let $U$ be an open subscheme of $S^0=S\setminus S_0$ and $p$ a prime invertible on $U$.
Fix an integer $m>0$.

Let $G$ be a semi-abelian scheme over $S$. Namely, $G$ is a separated smooth commutative
group scheme over $S$ such that each fiber $G_s$ of $G$ at $s\in S$ is an extension of
an abelian variety $A_s$ by a torus $T_s$. 
We denote the relative dimension of $G$ over $S$ by $d$.
Assume the following:

\begin{itemize}
 \item The rank of $T_s$ (called the toric rank of $G_s$) with $s\in S_0$ is a constant $r$.
 \item $G_U=G\times_SU$ is an abelian scheme.
\end{itemize}

Under the first condition, it is known that $G_0$ is globally an extension
\[
 0\longrightarrow T_0\longrightarrow G_0\longrightarrow A_0\longrightarrow 0,
\]
where $T_0$ is a torus of rank $r$ over $S_0$ and $A_0$ is an abelian scheme over $S_0$
(\cite[Chapter I, Corollary 2.11]{MR1083353}).

Let us consider two group spaces $\widehat{G}^\rig[p^m]_{U^\ad}$ and $G^\ad[p^m]_{U^\ad}$ over $U^\ad$, where $(-)_{U^\ad}$ denotes the restriction to $U^\ad$. 

\begin{lem}\label{lem:ad-fin-etale}
 The adic space $G^\ad[p^m]_{U^\ad}$ is finite \'etale of degree $p^{2dm}$ over $U^\ad$.
\end{lem}

\begin{prf}
 By Lemma \ref{lem:fiber-product}, we have $G^\ad[p^m]_{U^\ad}=(G_U[p^m])\times_UU^\ad$.
 Since $G_U[p^m]$ is finite \'etale of degree over $p^{2dm}$ over $U$, 
 $G^\ad[p^m]_{U^\ad}$ is finite \'etale of degree $p^{2dm}$ over $U^\ad$
 (\cf \cite[Corollary 1.7.3 i)]{MR1734903}).
\end{prf}

\begin{lem}\label{lem:rig-fin-etale}
 The adic space $\widehat{G}^\rig[p^m]_{U^\ad}$ is finite \'etale
 of degree $p^{(2d-r)m}$ over $U^\ad$.
\end{lem}

\begin{prf}
 We may assume that $S=\Spec R$ is affine. Let $I\subset R$ be the defining ideal of $S_0$.
 By replacing $R$ by its $I$-adic completion, we can reduce to the case where $R$ is $I$-adic complete.
 Put $S_i=\Spec R/I^{i+1}$ and $G_i=G\times_SS_i$.

 By \cite[Expos\'e IX, Th\'eor\`eme 3.6, Th\'eor\`eme 3.6 bis]{SGA3}, the exact sequence
 \[
  0\longrightarrow T_0\longrightarrow G_0\longrightarrow A_0\longrightarrow 0
 \]
 can be lifted canonically to an exact sequence
 \[
  0\longrightarrow T_i\longrightarrow G_i\longrightarrow A_i\longrightarrow 0
 \]
 over $S_i$, where $T_i$ is a torus over $S_i$ and $A_i$ is an abelian scheme over $S_i$
 (\cf \cite[\S 3.3.3]{Kai-Wen}). 
 Let $\widehat{T}=\varinjlim_i T_i$ and 
 $\widehat{A}=\varinjlim_i A_i$ be associated formal groups over $\mathcal{S}$.
 Then $\widehat{G}$ is an extension of $\widehat{A}$ by $\widehat{T}$.

 By taking $p^m$-torsion points, we get an exact sequence
 \[
  0\longrightarrow \widehat{T}[p^m]\longrightarrow \widehat{G}[p^m]\longrightarrow \widehat{A}[p^m]\longrightarrow 0
 \]
 of formal groups over $\mathcal{S}$. 
 Since $\widehat{G}^\rig[p^m]\cong (\widehat{G}[p^m])^\rig$, it suffices to see that 
 $(\widehat{T}[p^m])^\rig_{U^\ad}$ (resp.\ $(\widehat{A}[p^m])^\rig_{U^\ad}$) is finite \'etale of
 degree $p^{rm}$ (resp.\ $p^{2(d-r)m}$) over $U^\ad$.

 First we consider $(\widehat{T}[p^m])^\rig_{U^\ad}$. 
 Since $\widehat{T}[p^m]=\varinjlim_i (T_i[p^m])$, it is finite flat over $\mathcal{S}=\Spf R$.
 Therefore there exists a finite flat $R$-algebra $R'$ such that $\widehat{T}[p^m]=\Spf R'$.
 Moreover, a scheme $T'=\Spec R'$ is naturally equipped with a structure of commutative group scheme
 over $S=\Spec R$. Since $T'$ is killed by $p^m$ and $p$ is invertible on $U$, $T'_U=T'\times_SU$
 is a finite \'etale group scheme over $U$. 
 By Lemma \ref{lem:rig-ad} ii), we have
 $(\widehat{T}[p^m])^\rig=(T')^{\wedge\rig}=T'^\ad=T'\times_S\mathcal{S}^\rig$.
 Therefore $(\widehat{T}[p^m])^\rig_{U^\ad}=T'_U\times_UU^\ad$ is finite \'etale over $U^\ad$
 (\cf \cite[Corollary 1.7.3 i)]{MR1734903}). Its degree is clearly $p^{rm}$.

 The same argument also works for $(\widehat{A}[p^m])^\rig_{U^\ad}$.
\end{prf}

By Lemma \ref{lem:ad-fin-etale} and Lemma \ref{lem:rig-fin-etale}, we may regard
$\widehat{G}^\rig[p^m]_{U^\ad}$ and $G^\ad[p^m]_{U^\ad}$ as locally constant constructible sheaves
over $U^\ad$. 
Since we have a natural open immersion $\widehat{G}^\rig\hooklongrightarrow G^\ad$,
$\widehat{G}^\rig[p^m]_{U^\ad}$ is a subsheaf of $G^\ad[p^m]_{U^\ad}$.

\begin{rem}\label{rem:etale-sheaf-base-change}
 In the setting of Remark \ref{rem:rig-ad-base-change},
 the construction above is clearly compatible with the base change by $f\colon S'\longrightarrow S$.
\end{rem}

In the remaining part of this subsection, consider the case where $S=\Spec R$ is the spectrum of
a complete discrete valuation ring $R$,
$S_0$ is the closed point of $S$ and $U=S^0=S\setminus S_0$. 
Let $\overline{\eta}$ be a geometric point lying over the unique point of $U^{\ad}$.

As in the proof of Lemma \ref{lem:rig-fin-etale}, $\widehat{G}$ is an extension
\[
 0\longrightarrow \widehat{T}\longrightarrow \widehat{G}\longrightarrow \widehat{A}\longrightarrow 0
\]
of a formal group $\widehat{A}$ by $\widehat{T}$. 

\begin{lem}\label{lem:direct-summand}
 The $\Z/p^m\Z$-modules $\widehat{G}^\rig[p^m]_{\overline{\eta}}$ and 
 $G^\ad[p^m]_{\overline{\eta}}$ are free of finite rank. Furthermore, 
 $\widehat{G}^\rig[p^m]_{\overline{\eta}}$ (resp.\ $T_p\widehat{G}^\rig_{\overline{\eta}}$)
 is a direct summand of $G^\ad[p^m]_{\overline{\eta}}$
 (resp.\ $T_pG^\ad_{\overline{\eta}}$) as a $\Z_p$-module.
\end{lem}

\begin{prf}
 First $G^\ad[p^m]_{\overline{\eta}}\cong G_{\overline{\eta}}[p^m]$ is a free $\Z/p^m\Z$-module of finite rank,
 since $G_{\overline{\eta}}$ is an abelian variety. 

 Next we will prove that $\widehat{G}^\rig[p^m]_{\overline{\eta}}$ is a free $\Z/p^m\Z$-module of finite rank.
 By the exact sequence
 \[
  0\longrightarrow \widehat{T}^\rig[p^m]_{\overline{\eta}}\longrightarrow \widehat{G}^\rig[p^m]_{\overline{\eta}}\longrightarrow \widehat{A}^\rig[p^m]_{\overline{\eta}}\longrightarrow 0,
 \]
 it suffices to show that 
 $\widehat{T}^\rig[p^m]_{\overline{\eta}}$ and $\widehat{A}^\rig[p^m]_{\overline{\eta}}$ are free
 of finite rank.
 By \cite[Expos\'e X, Th\'eor\`eme 3.2]{SGA3}, $\widehat{T}$ can be algebraized into a torus $T$ over $S$.
 Therefore, by Lemma \ref{lem:fiber-product} and Lemma \ref{lem:rig-ad} ii), we have
 \[
  \widehat{T}^\rig[p^m]_{\overline{\eta}}=(T[p^m])^{\wedge\rig}_{\overline{\eta}}=(T[p^m])^\ad_{\overline{\eta}}
 =T_{\overline{\eta}}[p^m],
 \]
 which is obviously free of finite rank.
 Similar argument works for $\widehat{A}^\rig[p^m]_{\overline{\eta}}$, since 
 $\widehat{A}$ is also algebraizable (\cf \cite[Proposition 3.3.3.6, Remark 3.3.3.9]{Kai-Wen}). 
 
 In particular, $\widehat{G}^\rig[p^m]_{\overline{\eta}}$ is an injective $\Z/p^m\Z$-module,
 and thus is a direct summand of $G^\ad[p^m]_{\overline{\eta}}$.

 Put $C=T_pG^\ad_{\overline{\eta}}/T_p\widehat{G}^\rig_{\overline{\eta}}$.
 Then it is easy to see that 
 \[
  C\otimes_{\Z_p}\Z/p^m\Z=G^\ad[p^m]_{\overline{\eta}}/\widehat{G}^\rig[p^m]_{\overline{\eta}}
 \]
 (use the same argument as above to conclude that
 $T_p\widehat{G}^\rig_{\overline{\eta}}\otimes_{\Z_p}\Z/p^m\Z\cong \widehat{G}^\rig[p^m]_{\overline{\eta}}$).
 Therefore $C\otimes_{\Z_p}\Z/p^m\Z$ is a finite flat $\Z/p^m\Z$-module, and thus
 $C$ is a finite flat $\Z_p$-module. Hence $C$ is a free $\Z_p$-module, and thus
 the inclusion $T_p\widehat{G}^\rig_{\overline{\eta}}\hooklongrightarrow T_pG^\ad_{\overline{\eta}}$ splits.
 Namely, $T_p\widehat{G}^\rig_{\overline{\eta}}$ is a direct summand of $T_pG^\ad_{\overline{\eta}}$, as
 desired.
\end{prf}

\begin{rem}\label{rem:rig-Raynaud}
 Actually, the extension $0\longrightarrow \widehat{T}\longrightarrow \widehat{G}\longrightarrow \widehat{A}\longrightarrow 0$ can be algebraized; namely, there exists an exact sequence
\[
 0\longrightarrow T\longrightarrow G^\natural\longrightarrow A\longrightarrow 0
\] 
 of commutative group schemes over $S$, where $T$ is a torus over $S$ and $A$ is an abelian scheme over $S$,
 such that its formal completion along the special fiber is isomorphic to the extension above
 (\cf \cite[Proposition 3.3.3.6, Remark 3.3.3.9]{Kai-Wen}). Such an extension is called the Raynaud extension
 associated to $G$.

 Our construction above is related to the Raynaud extension in the following way.
 First, we have a natural isomorphism
 $\widehat{G}^\rig[p^m]_{\overline{\eta}}\yrightarrow{\cong}(G^\natural)^\ad[p^m]_{\overline{\eta}}$,
 which is induced from an open immersion
 $\widehat{G}^\rig\cong (\widehat{G^\natural})^\rig\hooklongrightarrow (G^\natural)^\ad$ 
 (\cf Lemma \ref{lem:rig-ad} i)). Moreover, the image of $\widehat{G}^\rig[p^m]_{\overline{\eta}}\hooklongrightarrow G^\ad[p^m]_{\overline{\eta}}$ coincides with the image of the map
 $G^\natural[p^m]_{\overline{\eta}}\longrightarrow G[p^m]_{\overline{\eta}}$ in \cite[Corollary 4.5.3.12]{Kai-Wen}.
 Therefore, the latter statement of Lemma \ref{lem:direct-summand} also follows from \cite[Corollary 4.5.3.12]{Kai-Wen}.
\end{rem}

\begin{prop}\label{prop:orthogonality}
 Assume that the fraction field $F$ of $R$ is a finite extension of $\Q_p$. 
 Let $\lambda$ be a polarization of $G_U$. Then, alternating bilinear pairings
 \[
 \langle\ ,\ \rangle_\lambda\colon G[p^m]_{\overline{\eta}}\times G[p^m]_{\overline{\eta}}\longrightarrow \mu_{p^m},\quad
 \langle\ ,\ \rangle_\lambda\colon T_pG_{\overline{\eta}}\times T_pG_{\overline{\eta}}\longrightarrow \Z_p(1)
 \]
 are induced. Then $\widehat{G}^\rig[p^m]_{\overline{\eta}}$
 (resp.\ $T_p\widehat{G}^\rig_{\overline{\eta}}$) is a coisotropic direct summand of 
 $G^\ad[p^m]_{\overline{\eta}}$ (resp.\ $T_pG^\ad_{\overline{\eta}}$) with respect to
 $\langle\ ,\ \rangle_\lambda$. Namely, we have
 $\widehat{G}^\rig[p^m]_{\overline{\eta}}^\perp\subset \widehat{G}^\rig[p^m]_{\overline{\eta}}$
 (resp.\ $(T_p\widehat{G}^\rig_{\overline{\eta}})^\perp\subset T_p\widehat{G}^\rig_{\overline{\eta}}$).
\end{prop}

\begin{prf}
 Put $V_p(-)=T_p(-)\otimes_{\Z_p}\Q_p$. Then an alternating bilinear pairing
 \[
  \langle\ ,\ \rangle_\lambda\colon V_pG_{\overline{\eta}}\times V_pG_{\overline{\eta}}\longrightarrow \Q_p(1)
 \]
 is induced. It suffices to show that $V_p\widehat{G}^\rig_{\overline{\eta}}$
 is a coisotropic subspace of $V_p\widehat{G}^{\ad}_{\overline{\eta}}$.
 
 We will prove a more precise result 
 $(V_p\widehat{G}^\rig_{\overline{\eta}})^\perp=V_p\widehat{T}^{\rig}_{\overline{\eta}}$.
 Since 
 \[
  \dim_{\Q_p}V_p\widehat{T}^\rig_{\overline{\eta}}+\dim_{\Q_p} V_p\widehat{G}^\rig_{\overline{\eta}}
 =r+(2d-r)=2d=\dim_{\Q_p} V_pG^{\ad}_{\overline{\eta}}
 \]
 (\cf Lemma \ref{lem:rig-fin-etale} and its proof), it is sufficient to prove that
 $V_p\widehat{T}^\rig_{\overline{\eta}}\subset (V_p\widehat{G}^\rig_{\overline{\eta}})^\perp$.
 Namely, we should prove that the homomorphism 
 $V_p\widehat{T}^\rig_{\overline{\eta}}\otimes_{\Q_p}V_p\widehat{G}^\rig_{\overline{\eta}}\longrightarrow \Q_p(1)$ induced by $\langle\ ,\ \rangle_\lambda$ is zero.
 
 Since $V_pG_{\overline{\eta}}$ is a semi-stable representation of
 $\Gal(\overline{F}/F)$, so are $V_p\widehat{T}^\rig_{\overline{\eta}}$ and $V_p\widehat{G}^\rig_{\overline{\eta}}$.
 Denote the residue field of $F$ by $\kappa_F$ and put $q=\#\kappa_F$.
 Consider the action of $\varphi^{[\kappa_F:\F_p]}$ on
 $D_{\mathrm{st}}(V_p\widehat{T}^\rig_{\overline{\eta}})$
 and $D_{\mathrm{st}}(V_p\widehat{A}^\rig_{\overline{\eta}})$.
 Let $T$ and $A$ be as in the proof of Lemma \ref{lem:direct-summand}.
 Then we have $V_p\widehat{T}^\rig_{\overline{\eta}}\cong V_pT_{\overline{\eta}}$
 and $V_p\widehat{A}^\rig_{\overline{\eta}}\cong V_pA_{\overline{\eta}}$.
 Therefore every eigenvalue of $\varphi^{[\kappa_F:\F_p]}$ on $D_{\mathrm{st}}(V_p\widehat{T}^\rig_{\overline{\eta}})$ is a Weil $q^{-2}$-number (for the definition of Weil numbers, see \cite[p.~471]{MR2276777}).
 Similarly, every eigenvalue of $\varphi^{[\kappa_F:\F_p]}$ on 
 $D_{\mathrm{st}}(V_p\widehat{A}^\rig_{\overline{\eta}})$ is a Weil $q^{-1}$-number;
 this is a consequence of the Weil conjecture for crystalline cohomology of abelian varieties.
 Therefore, every eigenvalue of $\varphi^{[\kappa_F:\F_p]}$ on 
 $D_{\mathrm{st}}(V_p\widehat{T}^\rig_{\overline{\eta}}\otimes_{\Q_p}V_p\widehat{G}^\rig_{\overline{\eta}})$
 is either a Weil $q^{-4}$-number or
 a Weil $q^{-3}$-number.
 On the other hand, every eigenvalue of $\varphi^{[\kappa_F:\F_p]}$ on $D_{\mathrm{st}}(\Q_p(1))$ is equal to
 $q^{-1}$, which is a Weil $q^{-2}$-number.
 Hence any $\varphi$-homomorphism 
 $D_{\mathrm{st}}(V_p\widehat{T}^\rig_{\overline{\eta}}\otimes_{\Q_p}V_p\widehat{G}^\rig_{\overline{\eta}})\longrightarrow D_{\mathrm{st}}(\Q_p(1))$ is zero.
 Since the functor $D_{\mathrm{st}}$ is fully faithful, 
 any $\Gal(\overline{F}/F)$-equivariant homomorphism
 \[
  V_p\widehat{T}^\rig_{\overline{\eta}}\otimes_{\Q_p}V_p\widehat{G}^\rig_{\overline{\eta}}\longrightarrow \Q_p(1)
 \]
 is zero. This completes the proof.
\end{prf}

\section{Shimura varieties of PEL type}\label{sec:Sh-var}
\subsection{Notation for Shimura varieties of PEL type}\label{subsec:notation-Sh-var}
In this paper, we are interested in Shimura varieties of PEL type considered in \cite[\S 5]{MR1124982}
(see also \cite[\S 1.4]{Kai-Wen}).
We recall it briefly. Fix a prime $p$. Consider a 6-tuple $(B,\mathcal{O}_B,*,V,L,\langle\ ,\ \rangle)$, where
\begin{itemize}
 \item $B$ is a finite-dimensional simple $\Q$-algebra such that $B\otimes_\Q\Q_p$ is a product of matrix algebras over unramified extensions of $\Q_p$,
 \item $\mathcal{O}_B$ is an order of $B$ whose $p$-adic completion is a maximal order of $B\otimes_\Q\Q_p$,
 \item $*$ is a positive involution of $B$ (namely, an involution such that $\Tr(bb^*)>0$ for every $b\in B^\times$) which preserves $\mathcal{O}_B$,
 \item $V$ is a non-zero finite $B$-module,
 \item $L$ is a $\Z$-lattice of $V$ preserved by $\mathcal{O}_B$, and
 \item $\langle\ ,\ \rangle\colon V\times V\longrightarrow \Q$ is a non-degenerate alternating
       $*$-Hermitian pairing with respect to $B$-action such that $\langle x,y\rangle\in\Z$ for every $x,y\in L$,
       and that $L_p=L\otimes_\Z\Z_p$ is a self-dual lattice of $V_p=V\otimes_\Q\Q_p$.
\end{itemize}
From $(B,*,V,\langle\ ,\ \rangle)$, we define a simple $\Q$-algebra $C=\End_B(V)$ with a unique involution $\#$
satisfying $\langle cv,w\rangle=\langle v,c^\# w\rangle$ for every $c\in C$ and $v,w\in V$.
Moreover we define an algebraic group $G$ over $\Q$ by
\[
 G(R)=\{g\in (C\otimes_\Q R)^\times\mid gg^\#\in R^\times\}
\]
for every $\Q$-algebra $R$. The condition $gg^\#\in R^\times$ is equivalent to the existence of $c(g)\in R^\times$ such that
$\langle gv,gw\rangle=c(g)\langle v,w\rangle$ for every $v,w\in V\otimes_\Q R$.
By the presence of the lattice $L$, $G$ can be naturally extended to a group scheme over $\Z$, which is also
denoted by the same symbol $G$. For any ring $R$, we put $G_R=G\otimes_\Z R$.

Consider an $\R$-algebra homomorphism $h\colon \C\longrightarrow C\otimes_\Q\R$ preserving involutions (on $\C$, we consider the complex conjugation)
such that the symmetric real-valued bilinear form $(v,w)\longmapsto \langle v,h(i)w\rangle$ on $V\otimes_\Q\R$ is positive definite.
Such a 7-tuple $(B,\mathcal{O}_B,*,V,L,\langle\ ,\ \rangle,h)$ is said to be an unramified integral PEL datum.
Note that the map $h$ induces a homomorphism $\Res_{\C/\R}\mathbb{G}_m\longrightarrow G_\R$ of algebraic groups over $\R$, which is also denoted by $h$.

Let $F$ be the center of $B$ and $F^+$ be the subfield of $F$ consisting of elements fixed by $*$.
The existence of $h$ tells us that $N=[F:F^+](\dim_FC)^{1/2}/2$ is an integer.
An unramified integral PEL datum falls into the following three types:
\begin{description}
 \item[type (A)] $[F:F^+]=2$.
 \item[type (C)] $[F:F^+]=1$ and $C\otimes_\Q\R$ is isomorphic to a product of $M_{2N}(\R)$.
 \item[type (D)] $[F:F^+]=1$ and $C\otimes_\Q\R$ is isomorphic to a product of $M_N(\mathbb{H})$.
\end{description}
If we are in the case of type (D), we will exclude the prime $p=2$.

Using $h\colon \C\longrightarrow C\otimes_{\Q}\R\hooklongrightarrow C\otimes_{\Q}\C$,
we can decompose the $B\otimes_\Q\C$-module
$V\otimes_\Q\C$ as $V\otimes_\Q\C=V_1\oplus V_2$, where $V_1$ (resp.\ $V_2$) is the subspace of $V\otimes_\Q\C$
on which $h(z)$ acts by $z$ (resp.\ $\overline{z}$) for every $z\in\C$.
We denote by $E$ the field of definition of the isomorphism class of the $B\otimes_\Q\C$-module $V_1$,
and call it the reflex field. It is a subfield of $\C$ which is finite over $\Q$.

In the sequel, we fix an unramified integral PEL datum $(B,\mathcal{O}_B,*,V,L,\langle\ ,\ \rangle,h)$.
For a compact open subgroup $K^p$ of $G(\widehat{\Z}^p)$, consider the functor $\Sh_{K^p}$ from the category of
$\mathcal{O}_{E,(p)}$-schemes to the category of sets, that associates $S$ to the set of isomorphism classes
of quadruples $(A,i,\lambda,\eta^p)$, where
\begin{itemize}
 \item $A$ is an abelian scheme over $S$,
 \item $\lambda\colon A\longrightarrow A^\vee$ is a prime-to-$p$ polarization,
 \item $i\colon \mathcal{O}_B\longrightarrow \End_S(A)$ is an algebra homomorphism
       such that $\lambda\circ i(b)=i(b^*)^\vee\circ \lambda$ for every $b\in\mathcal{O}_B$,
 \item $\eta^p$ is a level-$K^p$ structure of $(A,i,\lambda)$ of type 
       $(L\otimes_\Z\widehat{\Z}^p,\langle\ ,\ \rangle)$ in the sense of \cite[Definition 1.3.7.6]{Kai-Wen},
\end{itemize}
satisfying the equality of polynomials $\det_{\mathcal{O}_S}(b;\Lie A)=\det_E(b;V^1)$ in the sense of
\cite[\S 5]{MR1124982}. Recall that two quadruples $(A,i,\lambda,\eta^p)$ and $(A',i',\lambda',\eta'^p)$
are said to be isomorphic if there exists an isomorphism $f\colon A\longrightarrow A'$ of abelian schemes
such that 
\begin{itemize}
 \item $\lambda=f^\vee\circ \lambda'\circ f$,
 \item $f\circ i(b)=i'(b)\circ f$ for every $b\in\mathcal{O}_B$,
 \item and $f\circ \eta^p=\eta'^p$ in the sense of \cite[Definition 1.4.1.4]{Kai-Wen}.
\end{itemize}
If $K^p$ is neat (\cf \cite[Definition 1.4.1.8]{Kai-Wen}), the functor $\Sh_{K^p}$ is represented by a quasi-projective smooth
$\mathcal{O}_{E,(p)}$-scheme \cite[Corollary 7.2.3.9]{Kai-Wen}, which is also denoted by $\Sh_{K^p}$.
In this paper, we will call it a Shimura variety of PEL type.
The group $G(\A^{p,\infty})$ naturally acts on the tower of schemes $(\Sh_{K^p})_{K^p\subset G(\widehat{\Z}^p)}$
as Hecke correspondences.

\begin{rem}
 \begin{enumerate}
  \item Our definition of $\Sh_{K^p}$, due to \cite{Kai-Wen}, is slightly different from that in \cite{MR1124982},
	but they give the same moduli space. See \cite[Proposition 1.4.3.4]{Kai-Wen}.
	We prefer this definition because we can speak about the universal abelian scheme $\mathcal{A}$
	over $\Sh_{K^p}$, rather than its isogeny class.
  \item The word ``Shimura variety'' is an abuse of terminology. Indeed, in some cases
	$\Sh_{K^p}\otimes_{\mathcal{O}_{E,(p)}}E$ is a finite disjoint union of
	Shimura varieties in the usual sense.
 \end{enumerate}
\end{rem}

Assume moreover the following condition due to \cite[Condition 1.4.3.10]{Kai-Wen}:
\begin{quote}
  The action of $\mathcal{O}_B$ on $L$ extends to an action of some maximal order $\mathcal{O}'_B$ of $B$ containing $\mathcal{O}_B$.
\end{quote}
As explained in the beginning of \cite[\S 5.2.1]{Kai-Wen}, this restriction is rather harmless.
Under this condition, thanks to \cite[Theorem 6.4.1.1]{Kai-Wen}, $\Sh_{K^p}$ has a toroidal compactification $\Sh^{\mathrm{tor}}_{K^p}$.
Furthermore, since $K^p$ is assumed to be neat,
we can (and do) take $\Sh^{\mathrm{tor}}_{K^p}$ so that it is a projective scheme over $\mathcal{O}_{E,(p)}$ \cite[Proposition 7.3.1.4, Theorem 7.3.3.4]{Kai-Wen}.
On $\Sh^{\mathrm{tor}}_{K^p}$, there exists a semi-abelian scheme extending the universal abelian scheme
$\mathcal{A}$ on $\Sh_{K^p}$. We denote it by the same symbol $\mathcal{A}$.

In the following, to ease notation, we fix a neat compact open subgroup $K^p$ of $G(\widehat{\Z}^p)$ and omit the subscript $K^p$.
For an integer $r\ge 0$, let $\Sh^{\mathrm{tor}}_r$ (resp.\ $\Sh^{\mathrm{tor}}_{\le r}$, resp.\ 
$\Sh^{\mathrm{tor}}_{\ge r}$) be the subset of $\Sh^{\mathrm{tor}}$ consisting of
$x\in \Sh^{\mathrm{tor}}$ such that the toric rank of the semi-abelian variety $\mathcal{A}_x$ is equal to
(resp.\ less than or equal to, resp.\  greater than or equal to) $r$.
By \cite[Lemma 3.3.1.4]{Kai-Wen}, $\Sh^{\mathrm{tor}}_{\le r}$ is an open subset of $\Sh^{\mathrm{tor}}$.
Therefore $\Sh^{\mathrm{tor}}_{\ge r}$ (resp.\ $\Sh^{\mathrm{tor}}_r$) is a closed (resp.\ locally closed)
subset of $\Sh^{\mathrm{tor}}$. We regard them as reduced subschemes of $\Sh^{\mathrm{tor}}$.
Obviously $\Sh^{\mathrm{tor}}_0=\Sh^{\mathrm{tor}}_{\le 0}$ coincides with the Shimura variety $\Sh$.

Now choose a place $v$ of $E$ dividing $p$ and write $\mathcal{O}_v$ for the ring of integers of $E_v$. 
By taking base change from $\mathcal{O}_{E,(p)}$ to $\mathcal{O}_v$, the schemes $\Sh, \Sh^\mathrm{tor},\ldots$
can be regarded as schemes over $\mathcal{O}_v$. We denote them by the same symbols. 
Moreover, we denote the fibers of them at $v$ by putting the subscripts $v$. Namely,
we set $\Sh_v=\Sh\otimes_{\mathcal{O}_v}\kappa_v$, 
$\Sh^{\mathrm{tor}}_{v,r}=\Sh^{\mathrm{tor}}_r\otimes_{\mathcal{O}_v}\kappa_v$ and so on, where $\kappa_v$
denotes the residue field of $\mathcal{O}_v$. Similarly,
we denote the generic fibers of them by putting the subscripts $\eta$.

Since $\Sh^{\mathrm{tor}}$ is proper over $\mathcal{O}_v$, we may consider the specialization map
\[
 \spp\colon (\Sh^{\mathrm{tor}}_\eta)^\ad=(\Sh^{\mathrm{tor}})^\ad=(\Sh^{\mathrm{tor}})^{\wedge\rig}\longrightarrow \Sh^{\mathrm{tor}}_v
\]
introduced in the end of Section \ref{subsec:notation-adic}
(for the second equality, see Lemma \ref{lem:rig-ad} ii)).

\begin{defn}\label{defn:rigid-Sh}
 For an integer $r\ge 0$, we put
 \[
 \Sh_{\mathopen{]}\ge r\mathclose{[}}=\spp^{-1}(\Sh^\mathrm{tor}_{v,\ge r})^\circ\cap \Sh_\eta^\ad,\quad 
 \Sh_{\mathopen{]}r\mathclose{[}}=\Sh_{\mathopen{]}\ge r\mathclose{[}}\setminus \Sh_{\mathopen{]}\ge r+1\mathclose{[}},
 \]
 where $(-)^\circ$ denotes the interior in $(\Sh^{\mathrm{tor}})^\ad$.
 Moreover, we define an open subset $\Sh_{\mathopen{]}r\mathclose{[}}^\vartriangle$
 of $\Sh_{\mathopen{]}r\mathclose{[}}$ by
 \[
 \Sh_{\mathopen{]}r\mathclose{[}}^\vartriangle=\spp^{-1}(\Sh^\mathrm{tor}_{v,r})^\circ\cap \Sh_\eta^\ad=\Sh_{\mathopen{]}\ge r\mathclose{[}}\cap \spp^{-1}(\Sh^{\mathrm{tor}}_{v,\le r}).
 \]
\end{defn}

So far, we have only considered level structures which are prime to $p$. Now we add $p^m$-level structures
on the universal abelian scheme of the generic fiber $\Sh_\eta$.
Let $\Sh_{m,\eta}$ be the scheme over $\Sh_\eta$ classifying principal level-$m$ structures
(\cf \cite[Definition 1.3.6.2]{Kai-Wen})
of the universal object $(\mathcal{A},i^{\mathrm{univ}},\lambda^{\mathrm{univ}})$ over $\Sh_\eta$.
We denote the structure morphism $\Sh_{m,\eta}\longrightarrow \Sh_{\eta}$ by $p_m$, which is finite and \'etale.

Let $K_m$ be the compact open subgroup of $G(\Q_p)$ defined as the kernel of $G(\Z_p)\longrightarrow G(\Z/p^m\Z)$. 
Then $\Sh_{m,\eta}$ coincides with the Shimura variety over $E_v$ with $K_mK^p$-level structures.
In particular, for an element $g$ of $G(\Q_p)$ and integers $m,m'\ge 0$ satisfying $g^{-1}K_mg\subset K_{m'}$,
we can define the Hecke operator $[g]\colon \Sh_{m,\eta}\longrightarrow \Sh_{m',\eta}$,
which is finite and \'etale.
For later use, let us recall its description.
Let $k$ be an extension field of $E_v$ and $x$ a $k$-valued point of $\Sh_{m,\eta}$ corresponding to $[(A,i,\lambda,\eta^p,\eta_p)]$.
Then $[g](x)=[(A',i',\lambda',\eta'^p,\eta'_p)]$ is given as follows.

First we assume that $L_p\subset gL_p$. By the condition on $\eta_p$, it lifts to an isomorphism
$\widetilde{\eta}_p\colon L_p\yrightarrow{\cong} T_pA_{\overline{k}}$, where $\overline{k}$ is
a separable closure of $k$.\label{Hecke}
Take an integer $a\ge 0$ such that $p^agL_p\subset L_p$, and put 
$A'_{\overline{k}}=A_{\overline{k}}/\widetilde{\eta}_{p,a}(p^agL_p/p^aL_p)$, where
$\widetilde{\eta}_{p,a}\colon L_p/p^aL_p\yrightarrow{\cong}A_{\overline{k}}[p^a]$ is the isomorphism
induced by $\widetilde{\eta}_p$. It is easy to see that $A'_{\overline{k}}$ is independent of the choice
of $a$. Moreover, the condition $g^{-1}K_mg\subset K_{m'}(\subset K_0)$ tells us that $A'_{\overline{k}}$ is also
independent of the choice of $\widetilde{\eta}_p$ (note that any other choice is $\widetilde{\eta}_p\circ h$
with $h\in K_m$, and that $h$ preserves $gL_p$). 
In particular, it descends to an abelian variety $A'$ over $k$.
As $\widetilde{\eta}_{p,a}(p^agL_p/p^aL_p)$ is stable under $\mathcal{O}_B$,
$A'$ carries an $\mathcal{O}_B$-multiplication $i'$ induced from $i$.
A natural $p$-isogeny $A\longrightarrow A'$ transfers the pair $(\lambda,\eta^p)$ on $A$
into a pair $(\lambda',\eta'^p)$ on $A'$.
The definition of $\eta'_p$ is slightly more complicated. By the construction of $A'_{\overline{k}}$,
the composite $L_p\yrightarrow{\widetilde{\eta}_p} T_pA_{\overline{k}}\longrightarrow T_pA'_{\overline{k}}$
extends uniquely to an isomorphism $\widetilde{\eta}'_p\colon gL_p\longrightarrow T_pA'_{\overline{k}}$.
Taking the mod-$p^{m'}$ reduction of the composite $\widetilde{\eta}'_p\circ g\colon L_p\longrightarrow T_pA'_{\overline{k}}$, we obtain an isomorphism $\eta'_p\colon L_p/p^{m'}L_p\yrightarrow{\cong}A'_{\overline{k}}[p^{m'}]$.
By the condition $g^{-1}K_mg\subset K_{m'}$, it is independent of the choice of $\widetilde{\eta}_p$ and thus
descends to an isomorphism $\eta'_p\colon L_p/p^{m'}L_p\yrightarrow{\cong}A'[p^{m'}]$, which gives
a $p^{m'}$-level structure on $A'$.

For a general $g\in G(\Q_p)$, we can take an integer $b$ such that $L_p\subset p^bgL_p$.
Then we have $[g](x)=[p^bg](x)$, and thus are reduced to the case above.

\bigbreak
As in Definition \ref{defn:rigid-Sh}, we introduce several important subspaces of $\Sh_{m,\eta}^\ad$.

\begin{defn}
 We put
 \[
 \Sh_{m,\mathopen{]}\ge r\mathclose{[}}=p_m^{-1}(\Sh_{\mathopen{]}\ge r\mathclose{[}}),\quad
 \Sh_{m,\mathopen{]}r\mathclose{[}}=p_m^{-1}(\Sh_{\mathopen{]}r\mathclose{[}}),\quad
 \Sh_{m,\mathopen{]}r\mathclose{[}}^\vartriangle=p_m^{-1}(\Sh_{\mathopen{]}r\mathclose{[}}^\vartriangle).
 \]
 Here we write the same symbol $p_m$ for the finite \'etale morphism
 $\Sh_{m,\eta}^\ad\longrightarrow \Sh_\eta^\ad$ induced by $p_m$.
\end{defn}

\begin{prop}\label{prop:torus-rank-Hecke}
 Let $g$ be an element of $G(\Q_p)$. For integers $m,m'\ge 0$ such that $g^{-1}K_m g\subset K_{m'}$,
 the Hecke operator $[g]\colon \Sh_{m,\eta}^\ad\longrightarrow \Sh_{m',\eta}^\ad$ maps
 $\Sh_{m,\mathopen{]}r\mathclose{[}}$ (resp.\ $\Sh_{m,\mathopen{]}\ge r\mathclose{[}}$)
 into $\Sh_{m',\mathopen{]}r\mathclose{[}}$ (resp.\ $\Sh_{m',\mathopen{]}\ge r\mathclose{[}}$).
\end{prop}

A point $x$ of an adic space $X$ which is locally of finite type over $\Spa(E_v,\mathcal{O}_v)$ is said to be classical
if the residue field of $x$ is a finite extension of $E_v$. For a subset $Y$ of $X$,
we denote the set of classical points in $Y$ by $Y(\cl)$.

First we consider about classical points in $\Sh_{m,\mathopen{]}r\mathclose{[}}$.

\begin{lem}\label{lem:torus-rank-Hecke-cl-pt}
 Let $x$ be a classical point of $\Sh_{m,\mathopen{]}r\mathclose{[}}$.
 Then $[g](x)\in \Sh_{m',\eta}^\ad$ lies in $\Sh_{m',\mathopen{]}r\mathclose{[}}$.
\end{lem}

\begin{prf}
 First note that $\Sh_{m,\mathopen{]}r\mathclose{[}}(\cl)=p_m^{-1}(\spp^{-1}(\Sh_{v,r}^{\mathrm{tor}})\cap \Sh_\eta^\ad)(\cl)$.
 This is an immediate consequence of \cite[Lemma 4.6.4 (c)]{MR1360610}, whose analogue for adic spaces
 can be proved similarly.

 Let $F$ be the residue field of $x$ and $\mathcal{O}_F$ its ring of integers.
 Since $\Sh^{\mathrm{tor}}$ is proper over $\mathcal{O}_v$, we have a natural morphism of schemes
 $\iota_x\colon \Spec \mathcal{O}_F\longrightarrow \Sh^{\mathrm{tor}}$. 
 The image of the generic point $\eta$ (resp.\ closed point $s$) of $\Spec \mathcal{O}_F$ under $\iota_x$
 can be identified with $x$ (resp.\ $\spp(x)$). Let us denote by $\mathcal{G}$ the pull-back of $\mathcal{A}$ under
 $\iota_x$. Then, $\mathcal{G}$ is a semi-abelian scheme over $\mathcal{O}_F$ whose fiber at $\eta$ (resp.\ $s$)
 is the abelian variety $\mathcal{A}_x$ (resp.\ a semi-abelian variety with toric rank $r$).

 Put $x'=[g](x)$, which is an $F$-valued point of $\Sh_{m',\eta}^\ad$.
 By the similar construction, we get a semi-abelian scheme $\mathcal{G}'$ over $\mathcal{O}_F$
 whose generic fiber $\mathcal{G}'_\eta$ coincides with $\mathcal{A}_{x'}$.
 It suffices to show that the toric rank of $\mathcal{G}'_s$ is equal to $r$.

 By the description of the Hecke operator $[g]$, $\mathcal{G}'_\eta=\mathcal{A}_{x'}$ is isomorphic to
 the quotient of $\mathcal{G}_\eta=\mathcal{A}_x$ by some finite \'etale subgroup scheme $C_\eta$.
 Let $C$ be the closure of $C_\eta$ in $\mathcal{G}$. Then, \cite[Lemma 3.4.3.1]{Kai-Wen} says that
 the quotient $\mathcal{G}/C$ is represented by a semi-abelian scheme over $\mathcal{O}_F$.
 Moreover, by \cite[Chapter I, Proposition 2.7]{MR1083353}, the isomorphism 
 $(\mathcal{G}/C)_\eta\cong \mathcal{G}'_\eta$ can be extended uniquely to an isomorphism 
 $\mathcal{G}/C\cong \mathcal{G}'$.
 Therefore $\mathcal{G}'_s$ is a quotient of $\mathcal{G}_s$ by a finite flat group scheme $C_s$.
 In particular, we have $\rank_{\Z_p}T_p\mathcal{G}_s=\rank_{\Z_p}T_p\mathcal{G}'_s$ and
 $\dim \mathcal{G}_s=\dim \mathcal{G}'_s$.
 This implies that the toric rank of $\mathcal{G}_s$ is the same as that of $\mathcal{G}'_s$.
\end{prf}

To deduce Proposition \ref{prop:torus-rank-Hecke} from Lemma \ref{lem:torus-rank-Hecke-cl-pt},
we will use the following easy lemma:

\begin{lem}\label{lem:cl-pt-dense}
 Let $X$ be an adic space locally of finite type over $\Spa(E_v,\mathcal{O}_v)$.
 \begin{enumerate}
  \item For every non-empty open subset $U$ of $X$, $U(\cl)$ is non-empty.
  \item For every open subset $U$, $U^-(\cl)$ is dense in $U^-$. Here $U^-$ denotes the closure of $U$ in $X$.
 \end{enumerate}
\end{lem}

\begin{prf}
 For i), we may assume that $U=\Spa(A,A^\circ)$, where $A$ is an affinoid algebra over $E_v$ and
 $A^\circ$ is the subring of $A$ consisting of power-bounded elements. Since $A\neq 0$, 
 $A$ has a maximal ideal $\mathfrak{m}$. By \cite[Corollary 6.1.2/3]{MR746961}, $F=A/\mathfrak{m}$ is
 a finite extension of $E_v$. Extend the valuation of $E_v$ to $F$. Then
 $v\colon A\longrightarrow A/\mathfrak{m}=F\longrightarrow \R$ gives a classical point of
 $U=\Spa(A,A^\circ)$.

 ii) is an immediate consequence of i); indeed, $U(\cl)$ is dense in $U^-$.
\end{prf}

\begin{prf}[of Proposition \ref{prop:torus-rank-Hecke}]
 Denote the complement of $\Sh_{m,\mathopen{]}\ge r\mathclose{[}}$ in $\Sh_{m,\eta}^\ad$
 by $(\Sh_{m,\mathopen{]}\ge r\mathclose{[}})^c$. 
 It coincides with $p_m^{-1}(\spp^{-1}(\Sh_{v,\le r-1})^-\cap \Sh_\eta^\ad)$, which is
 the closure of the open subset $p_m^{-1}(\spp^{-1}(\Sh_{v,\le r-1})\cap \Sh_\eta^\ad)$ of $\Sh_{m,\eta}^\ad$
 (note that $p_m$ is an open map).
 Therefore, $(\Sh_{m,\mathopen{]}\ge r\mathclose{[}})^c(\cl)$ is dense in
 $(\Sh_{m,\mathopen{]}\ge r\mathclose{[}})^c$ by Lemma \ref{lem:cl-pt-dense}.
 Similarly, $([g]^{-1}(\Sh_{m',\mathopen{]}\ge r\mathclose{[}})^c)(\cl)$ is dense
 in $[g]^{-1}(\Sh_{m',\mathopen{]}\ge r\mathclose{[}})^c$.

 On the other hand, Lemma \ref{lem:torus-rank-Hecke-cl-pt} says that 
 $(\Sh_{m,\mathopen{]}\ge r\mathclose{[}})^c(\cl)=([g]^{-1}(\Sh_{m',\mathopen{]}\ge r\mathclose{[}})^c)(\cl)$.
 From it we can deduce $(\Sh_{m,\mathopen{]}\ge r\mathclose{[}})^c=[g]^{-1}(\Sh_{m',\mathopen{]}\ge r\mathclose{[}})^c$, as both sides are closed in $\Sh_{m,\eta}^\ad$.
 Thus we have $\Sh_{m,\mathopen{]}\ge r\mathclose{[}}=[g]^{-1}(\Sh_{m',\mathopen{]}\ge r\mathclose{[}})$ and
 $\Sh_{m,\mathopen{]}r\mathclose{[}}=[g]^{-1}(\Sh_{m',\mathopen{]}r\mathclose{[}})$.
\end{prf}

\subsection{The set $\mathcal{S}_m$}
This subsection provides some basic preliminaries for the next subsection.
We denote by $\mathcal{S}_\infty$ the set of $\mathcal{O}_B\otimes_\Z\Z_p$-submodules $I\subset L_p$ such that
\begin{itemize}
 \item $I$ is a non-trivial direct summand of $L_p$ as an $\mathcal{O}_B\otimes_\Z\Z_p$-module, and
 \item $I$ is coisotropic, namely, $I^\perp\subset I$.
\end{itemize}
For an element $I\in\mathcal{S}_\infty$ and an integer $m\ge 1$, we have an
$\mathcal{O}_B\otimes_\Z\Z_p$-submodule $I/p^mI\subset L_p/p^mL_p$.
We denote the set consisting of such submodules by $\mathcal{S}_m$.

\begin{prop}\label{prop:S_m-parabolic}
 \begin{enumerate}
  \item The set $\mathcal{S}_\infty$ can be identified with the set $\mathcal{S}_{\Q_p}$ consisting of non-trivial
	coisotropic $B\otimes_\Q\Q_p$-submodules of $V_p$.
	In particular, $\mathcal{S}_\infty$ is naturally equipped with an action of $G(\Q_p)$.
  \item	There exist finitely many proper parabolic subgroups $P_1,\ldots,P_k$
	of $G_{\Q_p}$ such that $\mathcal{S}_\infty$ is isomorphic to $\coprod_{i=1}^k G(\Q_p)/P_i(\Q_p)$
	as $G(\Q_p)$-sets.
  \item We have a natural $G(\Z_p)$-equivariant bijection 
	$K_m\backslash \mathcal{S}_\infty\yrightarrow{\cong}\mathcal{S}_m$.
  \item Let $g$ be an element of $G(\Q_p)$. For integers $m,m'\ge 1$ with $g^{-1}K_mg\subset K_{m'}$, 
	we have a unique map $g^{-1}\colon \mathcal{S}_m\longrightarrow \mathcal{S}_{m'}$ that makes
	the following diagram commutative:
	\[
	 \xymatrix{%
	\mathcal{S}_\infty\ar@{->>}[r]\ar[d]^-{g^{-1}}& \mathcal{S}_m\ar[d]^-{g^{-1}}\\
	\mathcal{S}_\infty\ar@{->>}[r]& \mathcal{S}_{m'}\lefteqn{.}
	}
	\]
 \end{enumerate}
\end{prop}

\begin{prf}[of Proposition \ref{prop:S_m-parabolic} i)]
 First let $I$ be a $B\otimes_{\Q}\Q_p$-submodule of $V_p$.
 As $\mathcal{O}_B\otimes_\Z\Z_p$ is a maximal order of $B\otimes_\Q\Q_p$,
 it is hereditary \cite[Theorem 21.4]{MR0393100}.
 Therefore the $\mathcal{O}_B\otimes_\Z\Z_p$-module $L_p/I\cap L_p$ is projective, since
 it is finitely generated and torsion-free as a $\Z_p$-module (\cf \cite[Corollary 21.5]{MR0393100}).
 Therefore the canonical surjection $L_p\longrightarrow L_p/I\cap L_p$ splits, and
 $I\cap L_p$ is a direct summand of $L_p$ as an $\mathcal{O}_B\otimes_\Z\Z_p$-module.

 By the observation above, we can define the following two maps 
 \[
  \mathcal{S}_{\Q_p}\longrightarrow \mathcal{S}_\infty;\ I\longmapsto I\cap L_p,\qquad
 \mathcal{S}_\infty\longrightarrow \mathcal{S}_{\Q_p};\ I\longmapsto I\otimes_{\Z_p}\Q_p.
 \]
 It is easy to see that these are the inverse of each other.
\end{prf}

\begin{rem}\label{rem:projective-module}
 By the proof above, we know that $L_p$ is a projective $\mathcal{O}_B\otimes_{\Z}\Z_p$-module.
\end{rem}

For proving the remaining part, we need some preparation; the proof of Proposition \ref{prop:S_m-parabolic}
is completed in the end of this subsection. 

First let us introduce more notation.
Let $\Lambda$ be one of $\Z/p^m\Z$, $\Z_p$ or $\Q_p$, and put $\mathcal{O}_\Lambda=\mathcal{O}_B\otimes_\Z\Lambda$.
For the center $F$ of $B$, we have canonical decompositions
$F\otimes_\Q\Q_p=\prod_{\tau}F_\tau$ and $\mathcal{O}_F\otimes_\Z\Z_p=\prod_{\tau}\mathcal{O}_\tau$,
where $\tau$ runs through places of $F$ lying over $p$.
These decompositions are extended to the decomposition
$\mathcal{O}_{\Lambda}=\prod_{\tau}\mathcal{O}_{\tau,\Lambda}$.

By the assumption, $\mathcal{O}_{\tau,\Lambda}$ is isomorphic to a matrix algebra 
$M_{d_\tau}(\mathcal{O}_\tau\otimes_{\Z_p}\Lambda)$ over $\mathcal{O}_\tau\otimes_{\Z_p}\Lambda$.
As in \cite[\S 1.1.3]{Kai-Wen}, let us denote by $M_{\tau,\Lambda}$ the $\mathcal{O}_{\tau,\Lambda}$-module
$(\mathcal{O}_\tau\otimes_{\Z_p}\Lambda)^{d_\tau}$.
By \cite[Lemma 1.1.3.4]{Kai-Wen}, every finitely generated projective $\mathcal{O}_\Lambda$-module $M$
is isomorphic to $\bigoplus_{\tau}M_{\tau,\Lambda}^{m_\tau}$, for uniquely determined integers $m_\tau\ge 0$.
We call $\underline{m}=(m_\tau)$ the multi-rank of $M$. Put $\lvert \underline{m}\rvert=\sum_{\tau}m_\tau$.

Similarly, we can define the multi-rank of a finitely generated projective 
$B\otimes_{\Q}\overline{\Q}_p$-module.
For a finitely generated projective $B\otimes_\Q{\Q_p}$-module $M$, 
the multi-rank of the $B\otimes_{\Q}\overline{\Q}_p$-module $M\otimes_{\Q_p}\overline{\Q}_p$
is not equal to the multi-rank of $M$ in general.

\begin{lem}\label{lem:stab-parab}
 For $I\in \mathcal{S}_{\Q_p}$, the stabilizer $P_I$ of $I\subset V_p$ in $G_{\Q_p}$ is
 a proper parabolic subgroup of $G_{\Q_p}$.
\end{lem}

\begin{prf}
 Let $\underline{m}$ be the multi-rank of the $B\otimes_\Q\overline{\Q}_p$-module
 $I\otimes_{\Q_p}\overline{\Q}_p$ and
 $X$ the moduli space classifying coisotropic $B\otimes_\Q\overline{\Q}_p$-submodules
 of $V\otimes_\Q\overline{\Q}_p$ with multi-rank $\underline{m}$.
 Since $X$ is a closed subscheme of a Grassmannian, it is projective over $\overline{\Q}_p$.

 Let us prove that $G(\overline{\Q}_p)$ acts on $X(\overline{\Q}_p)$ transitively.
 Take an arbitrary coisotropic $B\otimes_\Q\overline{\Q}_p$-submodule $I'$
 of $V\otimes_\Q\overline{\Q}_p$ with multi-rank $\underline{m}$.
 It suffices to find $g\in G(\overline{\Q}_p)$ such that $I'=g(I\otimes_{\Q_p}\overline{\Q}_p)$.
 Put $J=(I\otimes_{\Q_p}\overline{\Q}_p)^\perp$ and $J'=I'^\perp$.
 These are totally isotropic $B\otimes_{\Q}\overline{\Q}_p$-submodules of $V\otimes_\Q\overline{\Q}_p$
 with the same multi-rank.
 It is sufficient to find $g\in G(\overline{\Q}_p)$ such that $J'=gJ$.
 Fix a $B\otimes_{\Q_p}\overline{\Q}_p$-isomorphism $f\colon J\yrightarrow{\cong}J'$.
 By \cite[Lemma 1.2.4.4]{Kai-Wen}, we can find orthogonal decompositions
 \[
 V\otimes_\Q\overline{\Q}_p=(J\oplus J^\vee)\oplus (J\oplus J^\vee)^\perp
 =(J'\oplus J'^\vee)\oplus (J'\oplus J'^\vee)^\perp.
 \]
 Since $(J\oplus J^\vee)^\perp$ and $(J'\oplus J'^\vee)^\perp$ have the same multi-rank,
 by \cite[Corollary 1.2.3.10]{Kai-Wen}, there exists a $B\otimes_{\Q_p}\overline{\Q}_p$-isomorphism 
 $h\colon (J\oplus J^\vee)^\perp\yrightarrow{\cong}(J'\oplus J'^\vee)^\perp$ preserving the
 alternating pairings.
 Then the composite
 \[
 V\otimes_\Q\overline{\Q}_p=(J\oplus J^\vee)\oplus (J\oplus J^\vee)^\perp\yrightarrow{(f\oplus (f^\vee)^{-1})\oplus h}(J'\oplus J'^\vee)\oplus (J'\oplus J'^\vee)^\perp=V\otimes_\Q\overline{\Q}_p
 \]
 gives an element $g$ of $G(\overline{\Q}_p)$ satisfying $gJ=J'$.

 Now we know that $X_\mathrm{red}$ is a homogeneous space for $G_{\overline{\Q}_p}$ in the sense of
 \cite[2.3.1]{MR1642713}. Furthermore, we can easily see that the $G_{\overline{\Q}_p}$-space $X_{\mathrm{red}}$ satisfies
 the conditions in \cite[Corollary 5.5.4]{MR1642713}. Therefore, by \cite[Theorem 5.5.5]{MR1642713},
 $X_{\mathrm{red}}$ is a quotient of $G_{\overline{\Q}_p}$ by $P_{I,\overline{\Q}_p}=P_I\otimes_{\Q_p}\overline{\Q}_p$. 
 Hence  $G_{\overline{\Q}_p}/P_{I,\overline{\Q}_p}$ is projective over $\overline{\Q}_p$,
 that is, $P_{I,\overline{\Q}_p}$ is a parabolic subgroup of $G_{\overline{\Q}_p}$.
 Thus we can conclude that $P_I$ is a parabolic subgroup of $G_{\Q_p}$.

 Finally we will observe $P_I\neq G_{\Q_p}$. Put $J=(I\otimes_{\Q_p}\overline{\Q}_p)^\perp$
 and consider the orthogonal decomposition
 $V\otimes_\Q\overline{\Q}_p=(J\oplus J^\vee)\oplus (J\oplus J^\vee)^\perp$ as above.
 It is easy to see that $(J^\vee)^\perp$ gives a $\overline{\Q}_p$-valued point of $X$ distinct from $I\otimes_{\Q_p}\overline{\Q}_p$.
 Therefore $G_{\overline{\Q}_p}/P_{I,\overline{\Q}_p}=X_{\mathrm{red}}$ has more than one
 closed point, and thus $P_I\neq G_{\Q_p}$.
\end{prf}

\begin{lem}\label{lem:orthogonal-decomp-lift}
 Let $(M,\langle\ ,\ \rangle)$ be a finitely generated self-dual symplectic projective $\mathcal{O}_{\Z_p}$-module
 in the sense of \cite[Definition 1.1.4.7]{Kai-Wen}, and $N\subset M$ a totally isotropic direct summand
 as an $\mathcal{O}_{\Z_p}$-module.
 Put $M_m=M\otimes_{\Z_p}\Z/p^m\Z$ and $N_m=N\otimes_{\Z_p}\Z/p^m\Z$
 for simplicity.

 By \cite[Lemma 1.2.4.4]{Kai-Wen}, we can find an embedding $N_m^\vee\hooklongrightarrow M_m$
 with totally isotropic image such that $M_m$ is decomposed into the orthogonal direct sum:
 \[
 M_m=(N_m\oplus N_m^\vee)\oplus (N_m\oplus N_m^\vee)^\perp.
 \]
 Let us fix such an embedding. Then, we can extend it to 
 an embedding $N^\vee\hooklongrightarrow M$ with totally isotropic image, 
 such that $M$ is decomposed into the orthogonal direct sum:
 \[
 M=(N\oplus N^\vee)\oplus (N\oplus N^\vee)^\perp.
 \]
\end{lem}

\begin{prf}
 We use the similar induction argument as in the proof of \cite[Proposition 1.2.4.6]{Kai-Wen}.
 If $N\neq 0$, we can find a direct sum decomposition $N=N'\oplus N''$ such that 
 the multi-rank $\underline{n}'$ of $N'$
 satisfies $\lvert \underline{n}'\rvert=1$
 (in the proof of \cite[Proposition 1.2.4.6]{Kai-Wen},
 such a decomposition is written as $M_1\cong M_{1,0}\oplus M_1'$).
 Let $\tau$ be a unique place of $F$ lying over $p$ with $n'_\tau=1$.
 Then $N'\cong M_{\tau,\Z_p}$ is generated by one element $x$.
 Moreover, $N'^\vee\cong M_{\tau\circ *,\Z_p}$ is also generated by a single element.
 Take a generator $\overline{y}$ of $N_m'^\vee$
 and regard it as an element of $M_m$ by the fixed embedding $N_m'^\vee\subset N_m^\vee\hooklongrightarrow M_m$.

 Let $(\kern-2.5pt|\ ,\ |\kern-2.5pt)\colon M\times M\longrightarrow \mathcal{O}_{\Z_p}$ be the skew-Hermitian
 pairing associated to $\langle\ ,\ \rangle$ (\cf \cite[Lemma 1.1.4.5]{Kai-Wen}).
 We have $(\kern-2.5pt|x,x|\kern-2.5pt)=0$, since $N$ is totally isotropic.
 Choose a lift $y\in M$ of $\overline{y}$. Then we have $(\kern-2.5pt|y,y|\kern-2.5pt)\in p^m\mathcal{O}_{\Z_p}$,
 since $N_m^\vee$ is totally isotropic. 
 As in the proof of \cite[Proposition 1.2.4.6]{Kai-Wen},
 we may assume that $(\kern-2.5pt|x,y|\kern-2.5pt)=1_{\tau\circ *}\in \mathcal{O}_{\tau\circ *,\Z_p}$.
 Moreover, replacing $y$ by $1_{\tau\circ *}y$, we may assume that $y=1_{\tau\circ *}y$.

 If $\tau\neq \tau\circ *$, put $y'=y$. 
 Then $(\kern-2.5pt|y',y'|\kern-2.5pt)=1_{\tau\circ *}(\kern-2.5pt|y',y'|\kern-2.5pt)1_\tau=0$
 and $(\kern-2.5pt|x,y'|\kern-2.5pt)=1_{\tau\circ *}$.
 Consider the case $\tau=\tau\circ *$.
 Since $(\kern-2.5pt|y,y|\kern-2.5pt)\in p^m\mathcal{O}_{\Z_p}$, we can find $b\in p^m\mathcal{O}_{\Z_p}$
 such that $b-b^*=(\kern-2.5pt|y,y|\kern-2.5pt)$ by \cite[Proposition 1.2.2.1]{Kai-Wen}.
 Put $y'=y+bx$. Then we have $(\kern-2.5pt|y',y'|\kern-2.5pt)=0$ and
 $(\kern-2.5pt|x,y'|\kern-2.5pt)=1_{\tau\circ *}$.

 Denote by $N'^\sim$ the $\mathcal{O}_{\tau\circ *,\Z_p}$-span of $y'$ in $M$. 
 The pairing $\langle\ ,\ \rangle$
 gives a natural homomorphism $N'^\sim\longrightarrow N'^\vee$, whose mod-$p^m$ reduction is the identity
 (we regard $N'^\vee_m$ as a submodule of $M$ by the fixed embedding $N'^\vee_m\subset N_m^\vee\hooklongrightarrow M_m$). Therefore it is an isomorphism and gives an embedding $N'^\vee\cong N'^\sim\hooklongrightarrow M$
 extending the fixed embedding $N'^\vee_m\subset N_m^\vee\hooklongrightarrow M_m$.
 By \cite[Lemma 1.2.3.9]{Kai-Wen}, we have an orthogonal decomposition
 \[
  M=(N'\oplus N'^\vee)\oplus (N'\oplus N'^\vee)^\perp.
 \]
 Since $N\subset M$ is totally isotropic, the $N'^\vee$-component of every element of $N''$ is zero.
 Therefore, if we write $N'''$ for the image of $N''$ under the projection 
 $M\longrightarrow (N'\oplus N'^\vee)^\perp$, we have $N=N'\oplus N'''$.
 As the image of the embedding $N_m''^\vee\subset N_m^\vee\hooklongrightarrow M_m$ is contained in
 $(N'_m\oplus N_m'^\vee)^\perp$, $N'''_m$ coincides with $N''_m$.
 Hence we can apply the induction hypothesis to
 $N'''\subset (N'\oplus N'^\vee)^\perp$ and $N_m'''^\vee=N_m''^\vee\hooklongrightarrow (N'_m\oplus N_m'^\vee)^\perp$. This concludes the proof.
\end{prf}

\begin{lem}\label{lem:reduction-isom}
 Let $(M,\langle\ ,\ \rangle)$ and $(N,\langle\ ,\ \rangle)$ be finitely generated
 self-dual symplectic projective $\mathcal{O}_{\Z_p}$-modules. Assume that for an integer $m\ge 1$ we are given
 an $\mathcal{O}_{\Z/p^m\Z}$-isomorphism $h\colon M\otimes_{\Z_p}\Z/p^m\Z\yrightarrow{\cong}N\otimes_{\Z_p}\Z/p^m\Z$ which preserves the alternating pairings up to $(\Z/p^m\Z)^\times$-multiplication.
 Then there exists an $\mathcal{O}_{\Z_p}$-isomorphism $\widetilde{h}\colon M\yrightarrow{\cong}N$
 preserving the alternating pairings up to $\Z_p^\times$-multiplication, such that
 $\widetilde{h}\bmod p^m=h$.
\end{lem}

\begin{prf}
 For a $\Z_p$-algebra $R$, let $H(R)$ (resp.\ $X(R)$) be the group (resp.\ set)
 consisting of $\mathcal{O}_{\Z_p}\otimes_{\Z_p}R$-automorphisms of $N\otimes_{\Z_p}R$
 (resp.\ $\mathcal{O}_{\Z_p}\otimes_{\Z_p}R$-isomorphisms $M\otimes_{\Z_p}R\yrightarrow{\cong}N\otimes_{\Z_p}R$)
 preserving the alternating pairings up to $R^\times$-multiplication.
 The functor $R\longmapsto H(R)$ (resp.\ $R\longmapsto X(R)$) is obviously represented by
 a group scheme (resp.\ scheme) over $\Z_p$. The group scheme $H$ naturally acts on $X$.

 Since $M\otimes_{\Z_p}\Z/p^m\Z$ is isomorphic to $N\otimes_{\Z_p}\Z/p^m\Z$, $M$ and $N$ are isomorphic
 as $\mathcal{O}_{\Z_p}$-modules (\cf \cite[Lemma 1.1.3.1]{Kai-Wen}).
 Therefore, by \cite[Corollary 1.2.3.10]{Kai-Wen}, $X$ is an $H$-torsor (with respect to the \'etale topology
 on $\Spec \Z_p$). On the other hand, by \cite[Corollary 1.2.3.12]{Kai-Wen}, $H$ is smooth over $\Z_p$.
 Thus we can conclude that $X$ is smooth over $\Z_p$. 
 In particular, the natural map $X(\Z_p)\longrightarrow X(\Z/p^m\Z)$ is surjective.
 Hence we can find $\widetilde{h}\in X(\Z_p)$ which is mapped to $h\in X(\Z/p^m\Z)$ under the map above,
 as desired.
\end{prf}

\begin{prf}[of Proposition \ref{prop:S_m-parabolic} ii), iii), iv)]
 First we prove iii).
 Since the natural surjection $\mathcal{S}_\infty\longrightarrow \mathcal{S}_m$ is $G(\Z_p)$-equivariant
 and $K_m$ acts trivially on $\mathcal{S}_m$, a $G(\Z_p)$-equivariant surjection
 $K_m\backslash \mathcal{S}_\infty\longrightarrow \mathcal{S}_m$ is induced.
 We should prove its injectivity. Let $I$ and $I'$ be two elements of $\mathcal{S}_\infty$
 such that $I_m=I'_m$ (as in the proof of Lemma \ref{lem:orthogonal-decomp-lift}, we write
 $I_m$ for $I\otimes_{\Z_p}\Z/p^m\Z$). 
 Put $N=I^\perp$ and $N'=I'^\perp$. They are totally isotropic direct summands of $L_p$.
 It suffices to find an element $g\in K_m$ such that $gN=N'$.

 By \cite[Lemma 1.2.4.4]{Kai-Wen}, Remark \ref{rem:projective-module} and Lemma \ref{lem:orthogonal-decomp-lift},
 we can find embeddings $i\colon N^\vee\hooklongrightarrow L_p$ and $i'\colon N'^\vee\hooklongrightarrow L_p$
 with totally isotropic images such that $i\bmod p^m=i'\bmod p^m$.
 These give the following orthogonal decompositions:
 \[
  L_p=(N\oplus N^\vee)\oplus (N\oplus N^\vee)^\perp=(N'\oplus N'^\vee)\oplus (N'\oplus N'^\vee)^\perp.
 \]
 Since we have
 \begin{align*}
  (N\oplus N^\vee)^\perp\otimes_{\Z_p}\Z/p^m\Z&=(N_m\oplus N_m^\vee)^\perp\\
  &=(N'_m\oplus N_m'^\vee)^\perp=(N'\oplus N'^\vee)^\perp\otimes_{\Z_p}\Z/p^m\Z,
 \end{align*}
 Lemma \ref{lem:reduction-isom} enables us to find an $\mathcal{O}_{\Z_p}$-isomorphism
 $h\colon (N\oplus N^\vee)^\perp\yrightarrow{\cong} (N'\oplus N'^\vee)^\perp$ preserving
 the alternating pairings up to multiplication by $a\in \Z_p^\times$ such that $h\bmod p^m=\id$.
 On the other hand, it is easy to find an $\mathcal{O}_{\Z_p}$-isomorphism
 $f\colon N\yrightarrow{\cong}N'$ such that $f\bmod p^m=\id$
 (\cf proof of Lemma \ref{lem:reduction-isom}).
 The composite
 \[
  L_p=(N\oplus N^\vee)\oplus (N\oplus N^\vee)^\perp
 \yrightarrow{(f\oplus a(f^\vee)^{-1})\oplus h}(N'\oplus N'^\vee)\oplus (N'\oplus N'^\vee)^\perp=L_p
 \]
 gives an element of $K_m$ satisfying $gN=N'$.
 This conclude the proof of iii).

 Next we prove ii). By Proposition \ref{prop:S_m-parabolic} i), we identify the $G(\Q_p)$-set 
 $\mathcal{S}_\infty$ with $\mathcal{S}_{\Q_p}$.
 By iii), $K_1\backslash\mathcal{S}_\infty$ is a finite set.
 In particular, $G(\Q_p)\backslash\mathcal{S}_\infty=G(\Q_p)\backslash\mathcal{S}_{\Q_p}$ is also a finite set.
 Let $I_1,\ldots,I_k\in \mathcal{S}_{\Q_p}$ be a system of representatives of
 $G(\Q_p)\backslash \mathcal{S}_{\Q_p}$. Let $P_i$ be the stabilizer of $I_i\subset V_p$ in $G_{\Q_p}$.
 Then we have an obvious isomorphism of $G(\Q_p)$-sets $\coprod_{i=1}^kG(\Q_p)/P_i(\Q_p)\cong \mathcal{S}_{\Q_p}$.
 Since $P_i$ is a proper parabolic subgroup of $G_{\Q_p}$ by Lemma \ref{lem:stab-parab},
 this completes the proof of ii).
  
 The remaining part iv) is an immediate consequence of iii). Indeed, we can define 
 $g^{-1}\colon K_m\backslash \mathcal{S}_{\infty}\longrightarrow K_m\backslash \mathcal{S}_{\infty}$ by
 $K_mI\longmapsto K_{m'}g^{-1}I$, which is well-defined as $g^{-1}K_mg\subset K_{m'}$.
\end{prf}

\subsection{Partition of $\Sh_{m,\mathopen{]}r\mathclose{[}}$}
Fix integers $m,r\ge 1$.
In this subsection, we give a partition of $\Sh_{m,\mathopen{]}r\mathclose{[}}$ indexed by $\mathcal{S}_m$.
First, let us apply the construction introduced in Section \ref{subsec:etale-sheaves-semi-abelian} to
the case where $S=\Sh^\mathrm{tor}_{\le r}$, $S_0=\Sh^\mathrm{tor}_{v,r}$, 
$U=\Sh_\eta$ and
$G=\mathcal{A}\vert_{\Sh^\mathrm{tor}_{\le r}}$. By Lemma \ref{lem:rig-gen-fiber}, 
$U^\ad=U\times_S t(\mathcal{S})_a=\Sh_{\mathopen{]}r\mathclose{[}}^\vartriangle$ 
in this case. Therefore, we have two locally constant constructible sheaves
\[
 \widehat{\mathcal{A}}^\rig[p^m]_{\Sh_{\mathopen{]}r\mathclose{[}}^\vartriangle}\subset
 \mathcal{A}^\ad[p^m]_{\Sh_{\mathopen{]}r\mathclose{[}}^\vartriangle}.
\]
By taking pull-back under 
$p_m\colon \Sh_{m,\mathopen{]}r\mathclose{[}}^\vartriangle=p_m^{-1}(\Sh_{\mathopen{]}r\mathclose{[}}^\vartriangle)\longrightarrow \Sh_{\mathopen{]}r\mathclose{[}}^\vartriangle$ and using the universal level structure
$p_m^*\mathcal{A}[p^m]\cong L_p/p^mL_p$ over $\Sh_{m,\mathopen{]}r\mathclose{[}}^\vartriangle$, we get
a locally constant subsheaf of $L_p/p^mL_p$ corresponding to $\widehat{\mathcal{A}}^\rig[p^m]_{\Sh_{\mathopen{]}r\mathclose{[}}^\vartriangle}$,
 which we denote by $\mathcal{F}_{m,r}$.

\begin{lem}
 For every $x\in \Sh_{m,\mathopen{]}r\mathclose{[}}^\vartriangle$, the stalk $\mathcal{F}_{m,r,\overline{x}}$
 is an element of $\mathcal{S}_m$.
\end{lem}

\begin{prf}
 Since $\Sh_{m,\mathopen{]}r\mathclose{[}}^\vartriangle(\cl)$ is dense in 
 $\Sh_{m,\mathopen{]}r\mathclose{[}}^\vartriangle$ (Lemma \ref{lem:cl-pt-dense}),
 we may assume that $x$ is a classical point.
 
 By definition, the level structure $\mathcal{A}[p^m]_{\overline{x}}\cong L_p/p^mL_p$ lifts to an isomorphism
 $T_p\mathcal{A}_{\overline{x}}\cong L_p$ of $\mathcal{O}_B\otimes_{\Z}\Z_p$-modules which preserves
 alternating bilinear pairings up to multiplication by $\Z_p^\times$; here $T_p\mathcal{A}_{\overline{x}}$ is
 equipped with the alternating bilinear pairing $T_p\mathcal{A}_{\overline{x}}\times T_p\mathcal{A}_{\overline{x}}\longrightarrow \Z_p(1)$ induced from the universal polarization. 
 Therefore, it suffices to show that $T_p\widehat{\mathcal{A}}^\rig_{\overline{x}}$ is a non-trivial coisotropic
 direct summand of $T_p\mathcal{A}^\ad_{\overline{x}}$ as an $\mathcal{O}_B\otimes_\Z\Z_p$-module.
 First note that $T_p\widehat{\mathcal{A}}^\rig_{\overline{x}}$ is an $\mathcal{O}_B\otimes_{\Z}\Z_p$-submodule of $T_p\mathcal{A}^\ad_{\overline{x}}$,
 since the semi-abelian scheme $\mathcal{A}$ over $\Sh^{\mathrm{tor}}$ is endowed with an $\mathcal{O}_B$-multiplication,
 which is an extension of the universal $\mathcal{O}_B$-multiplication on $\mathcal{A}\vert_{\Sh}$.

 By Lemma \ref{lem:direct-summand}, $T_p\widehat{\mathcal{A}}^\rig_{\overline{x}}$ is a direct summand of
 $T_p\mathcal{A}^\ad_{\overline{x}}$ as a $\Z_p$-module. 
 As we assume that $\mathcal{O}_B\otimes_\Z\Z_p$ is
 a maximal order of $B\otimes_\Q\Q_p$, it is hereditary. Therefore the quotient
 $T_p\mathcal{A}^\ad_{\overline{x}}/T_p\widehat{\mathcal{A}}^\rig_{\overline{x}}$ is a projective
 $\mathcal{O}_B\otimes_\Z\Z_p$-module, and thus $T_p\widehat{\mathcal{A}}^\rig_{\overline{x}}$ is
 a direct summand of $T_p\mathcal{A}^\ad_{\overline{x}}$ as an $\mathcal{O}_B\otimes_\Z\Z_p$-module.

 By Lemma \ref{lem:rig-fin-etale}, we have 
 $\rank_{\Z_p}T_p\widehat{\mathcal{A}}^\rig_{\overline{x}}=\dim_{\Q_p}V_p-r$.
 Since $1\le r\le 1/2\dim_{\Q_p}V_p$, $T_p\widehat{\mathcal{A}}^\rig_{\overline{x}}$ is a non-trivial submodule
 of $T_p\mathcal{A}^\ad_{\overline{x}}$.
 Finally, Proposition \ref{prop:orthogonality} tells us that 
 $(T_p\widehat{\mathcal{A}}^{\rig}_{\overline{x}})^\perp\subset T_p\widehat{\mathcal{A}}^\rig_{\overline{x}}$.
\end{prf}

\begin{defn}
 For $I\in\mathcal{S}_m$, let $\Sh_{m,\mathopen{]}r\mathclose{[},I}^\vartriangle$ be the subset of
 $\Sh_{m,\mathopen{]}r\mathclose{[}}^\vartriangle$ consisting of $x\in \Sh_{m,\mathopen{]}r\mathclose{[}}^\vartriangle$
 such that $\mathcal{F}_{m,r,\overline{x}}=I$. 
 It is open and closed in $\Sh_{m,\mathopen{]}r\mathclose{[}}^\vartriangle$.
 Denote the closure of $\Sh_{m,\mathopen{]}r\mathclose{[},I}^\vartriangle$ in 
 $\Sh_{m,\mathopen{]}r\mathclose{[}}$ by $\Sh_{m,\mathopen{]}r\mathclose{[},I}$.
\end{defn}

\begin{prop}\label{prop:partition-I}
 We have $\Sh_{m,\mathopen{]}r\mathclose{[}}=\coprod_{I\in\mathcal{S}_m}\Sh_{m,\mathopen{]}r\mathclose{[},I}$
 as topological spaces.
\end{prop}

\begin{prf}
 First note that 
 \[
  \Sh_{m,\mathopen{]}r\mathclose{[}}^\vartriangle=\Sh_{m,\mathopen{]}\ge r\mathclose{[}}\cap p_m^{-1}\spp^{-1}(\Sh^{\mathrm{tor}}_{v,\le r}),\quad 
 \Sh_{m,\mathopen{]}r\mathclose{[}}=\Sh_{m,\mathopen{]}\ge r\mathclose{[}}\cap \bigl(p_m^{-1}\spp^{-1}(\Sh^{\mathrm{tor}}_{v,\le r})\bigr)^-.
 \]
 Since $p_m^{-1}\spp^{-1}(\Sh^{\mathrm{tor}}_{v,\le r})$ is a quasi-compact open subset of
 the spectral space $(\Sh^{\mathrm{tor}}_{m,\eta})^\ad$, its closure 
 $(p_m^{-1}\spp^{-1}(\Sh^{\mathrm{tor}}_{v,\le r}))^-$
 consists of specializations of points in $p_m^{-1}\spp^{-1}(\Sh^{\mathrm{tor}}_{v,\le r})$
 (\cf \cite[Corollary of Theorem 1]{MR0251026}).
 Thus every point in $\Sh_{m,\mathopen{]}r\mathclose{[}}$ is a specialization of a point in $\Sh_{m,\mathopen{]}r\mathclose{[}}^\vartriangle$. Similarly, every point in $\Sh_{m,\mathopen{]}r\mathclose{[},I}$ is 
 a specialization of a point in $\Sh_{m,\mathopen{]}r\mathclose{[},I}^\vartriangle$;
 note that $\Sh_{m,\mathopen{]}r\mathclose{[},I}^\vartriangle=\Sh_{m,\mathopen{]}r\mathclose{[},I}\cap p_m^{-1}\spp^{-1}(\Sh^{\mathrm{tor}}_{v,\le r})\hooklongrightarrow \Sh_{m,\mathopen{]}r\mathclose{[},I}$ is a quasi-compact open immersion between locally spectral spaces.

 Now let $x$ be a point of $\Sh_{m,\mathopen{]}r\mathclose{[}}$. Then there exists a generalization 
 $y\in \Sh_{m,\mathopen{]}r\mathclose{[}}^\vartriangle$ of $x$.
 Since $\Sh_{m,\mathopen{]}r\mathclose{[}}^\vartriangle$ is
 a disjoint union of $\Sh_{m,\mathopen{]}r\mathclose{[},I}^\vartriangle$, there exists $I\in \mathcal{S}_m$ with
 $y\in \Sh_{m,\mathopen{]}r\mathclose{[},I}^\vartriangle$. 
 Then $x$ lies in $\Sh_{m,\mathopen{]}r\mathclose{[},I}$.

 Next we will prove that $\Sh_{m,\mathopen{]}r\mathclose{[},I}$ and 
 $\Sh_{m,\mathopen{]}r\mathclose{[},I'}$ are disjoint for $I,I'\in\mathcal{S}_m$ with $I\neq I'$.
 Assume that we can take a point $x$ in $\Sh_{m,\mathopen{]}r\mathclose{[},I}\cap \Sh_{m,\mathopen{]}r\mathclose{[},I'}$.
 Then there exist points $y\in \Sh_{m,\mathopen{]}r\mathclose{[},I}^\vartriangle$ and
 $y'\in \Sh_{m,\mathopen{]}r\mathclose{[},I'}^\vartriangle$ both of which are generalizations of $x$.
 By \cite[Lemma 1.1.10 i)]{MR1734903}, $y$ is a specialization of $y'$ or $y'$ is a specialization of $y$.
 Since $\Sh_{m,\mathopen{]}r\mathclose{[},I}^\vartriangle$ is stable under specializations and generalizations
 in $\Sh_{m,\mathopen{]}r\mathclose{[}}^\vartriangle$, $y$ lies in 
 $\Sh_{m,\mathopen{]}r\mathclose{[},I}^\vartriangle\cap \Sh_{m,\mathopen{]}r\mathclose{[},I'}^\vartriangle$.
 This is contradiction.

 We have obtained $\Sh_{m,\mathopen{]}r\mathclose{[}}=\coprod_{I\in\mathcal{S}_m}\Sh_{m,\mathopen{]}r\mathclose{[},I}$
 as sets. Since $\Sh_{m,\mathopen{]}r\mathclose{[},I}$ is closed in $\Sh_{m,\mathopen{]}r\mathclose{[}}$ and
 $\mathcal{S}_m$ is a finite set, $\Sh_{m,\mathopen{]}r\mathclose{[},I}$ is also open in
 $\Sh_{m,\mathopen{]}r\mathclose{[}}$. Therefore we have an equality of topological spaces
 $\Sh_{m,\mathopen{]}r\mathclose{[}}=\coprod_{I\in\mathcal{S}_m}\Sh_{m,\mathopen{]}r\mathclose{[},I}$, 
 as desired.
\end{prf}

\begin{prop}\label{prop:partition-I-Hecke}
 Let $g$ be an element of $G(\Q_p)$. For integers $m,m'\ge 1$ such that $g^{-1}K_m g\subset K_{m'}$,
 the Hecke operator $[g]\colon \Sh_{m,\eta}^\ad\longrightarrow \Sh_{m',\eta}^\ad$ maps
 $\Sh_{m,\mathopen{]}r\mathclose{[},I}$ into $\Sh_{m',\mathopen{]}r\mathclose{[},g^{-1}I}$
 (for the definition of $g^{-1}I$, see Proposition \ref{prop:S_m-parabolic} iv)).
\end{prop}

\begin{prf}
 By Lemma \ref{lem:cl-pt-dense}, $\Sh_{m,\mathopen{]}r\mathclose{[},I}^\vartriangle(\cl)$ is
 dense in $\Sh_{m,\mathopen{]}r\mathclose{[},I}^\vartriangle$. Therefore
 $\Sh_{m,\mathopen{]}r\mathclose{[},I}^\vartriangle(\cl)$ is dense in $\Sh_{m,\mathopen{]}r\mathclose{[},I}$.
 By this fact and Proposition \ref{prop:torus-rank-Hecke}, it suffices to show that $[g]$ maps
 a classical point $x\in \Sh_{m,\mathopen{]}r\mathclose{[},I}^\vartriangle$ to a point in
 $\Sh_{m',\mathopen{]}r\mathclose{[},g^{-1}I}^\vartriangle$.
 Note that $\Sh_{m,\mathopen{]}r\mathclose{[}}^\vartriangle(\cl)=p_m^{-1}(\spp^{-1}(\Sh_{v,r}^{\mathrm{tor}})\cap \Sh_\eta^\ad)(\cl)=\Sh_{m,\mathopen{]}r\mathclose{[}}(\cl)$ (\cf \cite[Lemma 4.6.4 (c)]{MR1360610}).
 Hence $[g](x)$ lies in $\Sh_{m',\mathopen{]}r\mathclose{[}}^\vartriangle$ 
 by Proposition \ref{prop:torus-rank-Hecke}.

 We can reduce the problem to the case where $L_p\subset gL_p$.
 We use the same notation as in the proof of Lemma \ref{lem:torus-rank-Hecke-cl-pt}.
 Then we have the following commutative diagram of adic spaces over $\Spa(F,\mathcal{O}_F)$:
 \[
  \xymatrix{%
 \widehat{\mathcal{G}}^\rig[p^m]\ar@{^(->}[r]\ar[d]& \mathcal{G}^\ad[p^m]\ar[d]\\
 \widehat{\mathcal{G}}'^\rig[p^m]\ar@{^(->}[r]& \mathcal{G}'^\ad[p^m]\lefteqn{.}
 }
 \]
 This induces the commutative diagram
 \[
  \xymatrix{%
 T_p\widehat{\mathcal{G}}^\rig_{\overline{\eta}}\ar@{^(->}[r]\ar[d]& T_p\mathcal{G}^\ad_{\overline{\eta}}\ar[d]&
 L_p\ar[l]_-{\widetilde{\eta}_p}^-{\cong}\ar@{^(->}[d]\\
 T_p\widehat{\mathcal{G}}'^\rig_{\overline{\eta}}\ar@{^(->}[r]& T_p\mathcal{G}'^\ad_{\overline{\eta}}
 & gL_p\ar[l]_-{\widetilde{\eta}'_p}^-{\cong}& L_p\ar[l]_-{g}^-{\cong}\lefteqn{.}
 }
 \]
 Here $\widetilde{\eta}_p$ and $\widetilde{\eta}'_p$ are the same as in the description of the Hecke operator
 $[g]$ (page \pageref{Hecke}). By definition, $\mathcal{F}_{m',r,[g](x)}$ is the mod-$p^{m'}$ reduction
 of $g^{-1}\widetilde{\eta}_p'^{-1}(T_p\widehat{\mathcal{G}}'^\rig_{\overline{\eta}})$.
 Therefore, the diagram above tells us that $\mathcal{F}_{m',r,[g](x)}=g^{-1}\mathcal{F}_{m,r,x}=g^{-1}I$,
 as desired.
\end{prf}

\section{Cohomology of Shimura varieties}\label{sec:cohomology}
\subsection{Compactly supported cohomology and nearby cycle cohomology}\label{subsec:cohomology-notation}
As in Section \ref{subsec:notation-Sh-var}, let us fix a neat compact open subgroup $K^p$ of $G(\widehat{\Z}^p)$.
Fix a prime $\ell$ which is different from $p$.
Let $\xi$ be an algebraic representation of $G_{\Q_\ell}$ on a finite-dimensional 
$\overline{\Q}_\ell$-vector space. It naturally defines a smooth $\overline{\Q}_\ell$-sheaf $\mathcal{L}_\xi$
(or $\mathcal{L}_{\xi,m,K^p}$, if we need to indicate $m$ and $K^p$) on $\Sh_{m,\eta}$ 
(\cf \cite[III.2]{MR1876802}).
Moreover, $\mathcal{L_\xi}$ is equivariant with respect to the Hecke action.
For example, let $g$ be an element of $G(\Q_p)$ satisfying $g^{-1}K_mg\subset K_{m'}$ for integers $m,m'\ge 0$.
Then we have a natural isomorphism $[g]^*\mathcal{L}_{\xi,m',K^p}\yrightarrow{\cong}\mathcal{L}_{\xi,m,K^p}$
of smooth sheaves over $\Sh_{m,K^p,\eta}$.

Consider the compactly supported cohomology
\[
 H^i_c(\Sh_{\infty,K^p,\overline{\eta}},\mathcal{L}_\xi)=\varinjlim_{m}H^i_c(\Sh_{m,K^p,\eta}\otimes_{E_v}\overline{E}_v,\mathcal{L}_\xi).
\]
The group $G(\Q_p)\times \Gal(\overline{E}_v/E_v)$ naturally acts on it.
By this action, $H^i_c(\Sh_{\infty,K^p,\overline{\eta}},\mathcal{L}_\xi)$ becomes an admissible/continuous representation of
$G(\Q_p)\times \Gal(\overline{E}_v/E_v)$ in the sense of \cite[I.2]{MR1876802}.

On the other hand, we may also consider the nearby cycle cohomology defined as follows:
\[
 H^i_c(\Sh_{\infty,K^p,\overline{v}},R\psi\mathcal{L}_\xi)=\varinjlim_{m}H^i_c\bigl(\Sh_{K^p,\overline{v}},R\psi (p_{m*}\mathcal{L}_\xi)\bigr),
\]
where we put $\Sh_{K^p,\overline{v}}=\Sh_{K^p,v}\otimes_{\kappa_v}\overline{\kappa}_v$.
Obviously the group $\Gal(\overline{E}_v/E_v)$ acts on it. The following lemma gives an action of
$G(\Q_p)$ on $H^i_c(\Sh_{\infty,K^p,\overline{v}},R\psi\mathcal{L}_\xi)$.

\begin{lem}\label{lem:nearby-gp-action}
 We have a natural action of $G(\Q_p)$ on $H^i_c(\Sh_{\infty,K^p,\overline{v}},R\psi\mathcal{L}_\xi)$.
 By this action, $H^i_c(\Sh_{\infty,K^p,\overline{v}},R\psi\mathcal{L}_\xi)$ becomes an admissible/continuous
 $G(\Q_p)\times \Gal(\overline{E}_v/E_v)$-representation.
\end{lem}

\begin{prf}
 To ease notation, we omit the subscript $K^p$.

 As in \cite[\S 6]{MR2169874},
 we can construct a tower $(\Sh_m)_{m\ge 0}$ of schemes over $\mathcal{O}_v$ with finite transition maps
 such that $\Sh_m$ gives an integral model of $\Sh_{m,\eta}$ and $\Sh_0=\Sh$.
 In this situation, we have 
 \[
 H^i_c\bigl(\Sh_{\overline{v}},R\psi (p_{m*}\mathcal{L}_\xi)\bigr)=H^i_c(\Sh_{m,\overline{v}},R\psi\mathcal{L}_\xi),
 \]
 where $\Sh_{m,\overline{v}}=\Sh_m\otimes_{\mathcal{O}_v}\overline{\kappa}_v$.

 Put $G^+(\Q_p)=\{g\in G(\Q_p)\mid g^{-1}L_p\subset L_p\}$. For $g\in G^+(\Q_p)$, let
 $e(g)$ be the minimal non-negative integer such that $gL_p\subset p^{-e(g)}L_p$.
 Then we can construct a tower $(\Sh_{m,g})_{m\ge e(g)}$ of schemes over $\mathcal{O}_v$
 and two morphisms
 \[
 \pr\colon \Sh_{m,g}\longrightarrow \Sh_m,\quad [g]\colon \Sh_{m,g}\longrightarrow \Sh_{m-e(g)}
 \]
 which are compatible with the transition maps. It is known that these are proper morphisms,
 $\pr$ induces an isomorphism on the generic fibers,
 and $[g]$ induces the Hecke action of $g$ on the generic fibers
 \cite[Proposition 16, Proposition 17]{MR2169874}. In particular, we have a canonical cohomological
 correspondence (\cf \cite[Expos\'e III]{SGA5}, \cite{MR1431137})
 \[
  c_g\colon [g]_\eta^*\mathcal{L}_{\xi,m-e(g)}\yrightarrow{\cong} \pr_\eta^*\mathcal{L}_{\xi,m}=R\pr_\eta^!\mathcal{L}_{\xi,m}.
 \]
 Let
 \[
  R\psi(c_g)\colon [g]_{\overline{v}}^*R\psi\mathcal{L}_{\xi,m-e(g)}\longrightarrow R\pr_{\overline{v}}^!R\psi\mathcal{L}_{\xi,m}
 \]
 be the specialization of $c_g$ (\cf \cite[\S 1.5]{MR1431137}, \cite[\S 6]{RZ-GSp4}).
 Since $[g]_{\overline{v}}$ is proper, this induces a homomorphism
 \[
  H^i_c(\Sh_{m-e(g),\overline{v}},R\psi\mathcal{L}_\xi)
 \yrightarrow{H_c^i(R\psi(c_g))}H^i_c(\Sh_{m,\overline{v}},R\psi\mathcal{L}_\xi).
 \]
 Taking the inductive limit, we get
 \[
  \gamma_g\colon H^i_c(\Sh_{\infty,K^p,\overline{v}},R\psi\mathcal{L}_\xi)
 \longrightarrow H^i_c(\Sh_{\infty,K^p,\overline{v}},R\psi\mathcal{L}_\xi).
 \]
 From an obvious relation $c_{gg'}=c_g\circ g^*c_{g'}$ for $g,g'\in G^+(\Q_p)$,
 we deduce $\gamma_{gg'}=\gamma_{g}\circ \gamma_{g'}$ (\cf \cite[Corollary 6.3]{RZ-GSp4}).
 On the other hand, by \cite[Proposition 16 (3), Proposition 17 (3)]{MR2169874}, $\gamma_{p^{-1}}$ is
 the identity. Since $G(\Q_p)$ is generated by $G^+(\Q_p)$ and $p$ as a monoid, we can extend $\gamma_g$
 to the whole $G(\Q_p)$. By \cite[Proposition 16 (4)]{MR2169874}, the restriction of this action
 to $K_0=G(\Z_p)$ coincides with the inductive limit of the natural action of $K_0$ on
 $H^i_c(\Sh_{m,\overline{v}},R\psi \mathcal{L}_\xi)$.
 In particular, it is a smooth action. Furthermore, for integers $m'\ge m\ge 1$, we have
 \[
  H^i_c(\Sh_{m',\overline{v}},R\psi \mathcal{L}_\xi)^{K_m/K_{m'}}=H^i_c(\Sh_{m,\overline{v}},R\psi \mathcal{L}_\xi)
 \]
 since $K_m/K_{m'}$ is a $p$-group (\cf \cite[Proposition 2.5]{non-cusp}). Taking inductive limit, we obtain
 \[
  H^i_c(\Sh_{\infty,K^p,\overline{v}},R\psi\mathcal{L}_\xi)^{K_m}=H^i_c(\Sh_{m,\overline{v}},R\psi \mathcal{L}_\xi).
 \]
 This implies that $H^i_c(\Sh_{\infty,K^p,\overline{v}},R\psi\mathcal{L}_\xi)$ is
 an admissible/continuous representation of $G(\Q_p)\times\Gal(\overline{E}_v/E_v)$.
\end{prf}

Now we can state our main theorem in this article.

\begin{thm}\label{thm:main-thm}
 The kernel and cokernel of the canonical homomorphism
 \[
  H^i_c(\Sh_{\infty,K^p,\overline{v}},R\psi\mathcal{L}_\xi)\longrightarrow H^i_c(\Sh_{\infty,K^p,\overline{\eta}},\mathcal{L}_\xi)
 \]
 (\cf \cite[Expos\'e XIII, (2.1.7.3)]{SGA7})
 are non-cuspidal, namely, they have no supercuspidal subquotient of $G(\Q_p)$.
 In particular, for an irreducible supercuspidal representation $\pi$ of $G(\Q_p)$,
 we have an isomorphism 
 \[
  H^i_c(\Sh_{\infty,K^p,\overline{\eta}},\mathcal{L}_\xi)[\pi]\cong H^i_c(\Sh_{\infty,K^p,\overline{v}},R\psi\mathcal{L}_\xi)[\pi].
 \]
\end{thm}

\begin{rem}\label{rem:cpt-ordinary}
 Let $H^i(\Sh_{\infty,K^p,\overline{\eta}},\mathcal{L}_\xi)=\varinjlim_{m}H^i(\Sh_{m,K^p,\eta}\otimes_{E_v}\overline{E}_v,\mathcal{L}_\xi)$ be the ordinary cohomology of our Shimura variety.
 This is also an admissible/continuous representation of $G(\Q_p)\times \Gal(\overline{E}_v/E_v)$.
 By using the minimal compactification of $\Sh_{m,K^p,\eta}$ and its natural stratification
 (\cf \cite[\S 3.7]{MR1149032}),
 it is easy to see that the kernel and the cokernel of
 the canonical homomorphism
 \[
  H^i_c(\Sh_{\infty,K^p,\overline{\eta}},\mathcal{L}_\xi)\longrightarrow H^i(\Sh_{\infty,K^p,\overline{\eta}},\mathcal{L}_\xi)
 \]
 are non-cuspidal as $G(\Q_p)$-representations (in fact, we can use the similar argument
 as in the next subsection).
 Therefore, to prove Theorem \ref{thm:main-thm}, it suffices to 
 show that the kernel and the cokernel of the composite
 \[
  H^i_c(\Sh_{\infty,K^p,\overline{v}},R\psi\mathcal{L}_\xi)\longrightarrow H^i_c(\Sh_{\infty,K^p,\overline{\eta}},\mathcal{L}_\xi)\longrightarrow H^i(\Sh_{\infty,K^p,\overline{\eta}},\mathcal{L}_\xi)
 \]
 are non-cuspidal.
\end{rem}

\begin{rem}\label{rem:cpt-IH}
 Let $\mathit{IH}^i(\Sh_{\infty,K^p,\overline{\eta}},\mathcal{L}_\xi)=\varinjlim_{m}H^i(\Sh^{\mathrm{min}}_{m,K^p,\eta}\otimes_{E_v}\overline{E}_v,j_{!*}\mathcal{L}_\xi)$ be the intersection cohomology of our Shimura variety,
 where $j\colon \Sh_{m,K^p,\eta}\hooklongrightarrow \Sh_{m,K^p,\eta}^{\mathrm{min}}$ denotes the minimal
 compactification of $\Sh_{m,K^p,\eta}$.
 Then, as in the previous remark, it is easy to see that the kernel and cokernel of the canonical homomorphism
 \[
  H^i_c(\Sh_{\infty,K^p,\overline{\eta}},\mathcal{L}_\xi)\longrightarrow
 \mathit{IH}^i(\Sh_{\infty,K^p,\overline{\eta}},\mathcal{L}_\xi)
 \]
 is non-cuspidal. Therefore, by Theorem \ref{thm:main-thm},
 we have an isomorphism
 \[
  \mathit{IH}^i(\Sh_{\infty,K^p,\overline{\eta}},\mathcal{L}_\xi)[\pi]\cong H^i_c(\Sh_{\infty,K^p,\overline{v}},R\psi\mathcal{L}_\xi)[\pi]
 \]
 for an irreducible supercuspidal representation $\pi$ of $G(\Q_p)$.
\end{rem}

\begin{cor}
 Put
\[
 H^i_c(\Sh_{\infty,\overline{\eta}},\mathcal{L}_\xi)=\varinjlim_{K^p}H^i_c(\Sh_{\infty,K^p,\overline{\eta}},\mathcal{L}_\xi),\quad
 H^i_c(\Sh_{\infty,\overline{v}},R\psi\mathcal{L}_\xi)=\varinjlim_{K^p}H^i_c(\Sh_{\infty,K^p,\overline{v}},R\psi\mathcal{L}_\xi).
\]
 These are $G(\A^\infty)\times \Gal(\overline{E}_v/E_v)$-representations.

 Let $\Pi$ be an irreducible admissible representation of $G(\A^\infty)$ such that
 $\Pi_p$ is a supercuspidal representation of $G(\Q_p)$.
 Then, $\Pi$ does not appear as a subquotient of the kernel or the cokernel of the canonical
 homomorphism
 $H^i_c(\Sh_{\infty,\overline{v}},R\psi\mathcal{L}_\xi)\longrightarrow H^i_c(\Sh_{\infty,\overline{\eta}},\mathcal{L}_\xi)$. In particular, we have an isomorphism of $\Gal(\overline{E}_v/E_v)$-representations
 \[
 H^i_c(\Sh_{\infty,\overline{\eta}},\mathcal{L}_\xi)[\Pi]\cong H^i_c(\Sh_{\infty,\overline{v}},R\psi\mathcal{L}_\xi)[\Pi].
 \]
\end{cor}

\begin{prf}
 Take a neat compact open subgroup $K^p\subset G(\widehat{\Z}^p)$ such that
 $\Pi^{K^p}\neq 0$. If $\Pi$ appears as a subquotient of the kernel or the cokernel of 
 $H^i_c(\Sh_{\infty,\overline{v}},R\psi\mathcal{L}_\xi)\longrightarrow H^i_c(\Sh_{\infty,\overline{\eta}},\mathcal{L}_\xi)$, then $\Pi_p$ appears as a subquotient of the kernel or the cokernel of
 $H^i_c(\Sh_{\infty,K^p,\overline{v}},R\psi\mathcal{L}_\xi)\longrightarrow H^i_c(\Sh_{\infty,K^p,\overline{\eta}},\mathcal{L}_\xi)$. This contradicts to Theorem \ref{thm:main-thm}.
\end{prf}

\begin{rem}
 If an irreducible admissible representation $\Pi$ of $G(\A^\infty)$ is supercuspidal at a prime $p'\neq p$,
 we can prove the similar result as above. This case is due to Tetsushi Ito and the second author.
 Its proof is easier, since we may use minimal compactifications over $\Z_p$.
\end{rem}

\subsection{Proof of main theorem}
The starting point of our proof of Theorem \ref{thm:main-thm} is the following rigid-geometric interpretation
of the nearby cycle cohomology.

\begin{lem}\label{lem:nearby-rigid}
 There exist natural isomorphisms 
 \[
 H^i_c(\Sh_{\infty,K^p,\overline{v}},R\psi\mathcal{L}_\xi)\cong 
 \varinjlim_{m}H^i_c\bigl(\spp^{-1}(\Sh_{v,0}^{\mathrm{tor}})_{\overline{\eta}},p_{m*}\mathcal{L}^\ad_\xi\bigr)\cong
 \varinjlim_{m}H^i_{(\Sh_{m,\mathopen{]}\ge 1\mathclose{[}})^c_{\overline{\eta}}}(\Sh^\ad_{m,\overline{\eta}},\mathcal{L}^\ad_\xi).
 \]
Here the subscript $\overline{\eta}$ denotes the base change from $E_v$ to $\overline{E}_v$,
$(\Sh_{m,\mathopen{]}\ge 1\mathclose{[}})^c$ denotes the complement of $\Sh_{m,\mathopen{]}\ge 1\mathclose{[}}$
in $\Sh_{m,\eta}^\ad$, and $\mathcal{L}^\ad_\xi$ denotes the $\ell$-adic sheaf on $\Sh_{m,\eta}^\ad$ induced from
$\mathcal{L}_\xi$. As in Section \ref{subsec:notation-Sh-var}, we omit the subscript $K^p$.

 By Proposition \ref{prop:torus-rank-Hecke}, $\varinjlim_{m}H^i_{(\Sh_{m,\mathopen{]}\ge 1\mathclose{[}})^c_{\overline{\eta}}}(\Sh^\ad_{m,\overline{\eta}},\mathcal{L}^\ad_\xi)$ is naturally equipped with an action of $G(\Q_p)$.
 Under this action, the isomorphism
 \[
  H^i_c(\Sh_{\infty,K^p,\overline{v}},R\psi\mathcal{L}_\xi)\cong 
 \varinjlim_{m}H^i_{(\Sh_{m,\mathopen{]}\ge 1\mathclose{[}})^c_{\overline{\eta}}}(\Sh^\ad_{m,\overline{\eta}},\mathcal{L}^\ad_\xi)
 \]
 is $G(\Q_p)$-equivariant.
\end{lem}

\begin{rem}
 Before proving it, we should clarify the definition of the $\ell$-adic cohomology of rigid spaces
 appearing in the lemma. We define them simply by taking projective limit of the cohomology with
 torsion coefficients and tensoring $\overline{\Q}_\ell$. Note that this definition is, at least a priori,
 different from the compactly supported $\ell$-adic cohomology introduced in \cite{MR1626021}.
 However our naive definition also works well, since all the cohomology groups we encounter are known to be
 finitely generated.
\end{rem}

\begin{prf}[of Lemma \ref{lem:nearby-rigid}]
 Let $\mathcal{F}$ be an arbitrary $\ell^n$-torsion sheaf on $\Sh_{m,\overline{\eta}}$.
 We shall observe that there exist natural isomorphisms
 \[
  H^i_c(\Sh_{\overline{v}},R\psi p_{m*}\mathcal{F})\cong 
 H^i_c\bigl(\spp^{-1}(\Sh_{v,0}^{\mathrm{tor}})_{\overline{\eta}},p_{m*}\mathcal{F}^\ad\bigr)\cong
 H^i_{(\Sh_{m,\mathopen{]}\ge 1\mathclose{[}})^c_{\overline{\eta}}}(\Sh^\ad_{m,\overline{\eta}},\mathcal{F}^\ad).
 \]
 By \cite[Theorem 5.7.6, Theorem 3.7.2]{MR1734903}, we obtain the first isomorphism
 as the composite of
 \begin{align*}
  H^i_c(\Sh_{\overline{v}},R\psi p_{m*}\mathcal{F})&\cong H^i_c\bigl((\Sh^{\wedge\rig})_{\overline{\eta}},(p_{m*}\mathcal{F})^\ad\bigr)=H^i_c\bigl(\spp^{-1}(\Sh_{v,0}^{\mathrm{tor}})_{\overline{\eta}},p_{m*}\mathcal{F}^\ad\bigr).
 \end{align*}
 To construct the second isomorphism, we prove the following:
 \begin{itemize}
  \item $(\Sh_{\mathopen{]}\ge 1\mathclose{[}})^c$ contains $\spp^{-1}(\Sh_{v,0}^{\mathrm{tor}})$, and
  \item $(\Sh_{\mathopen{]}\ge 1\mathclose{[}})^c$ is closed in $(\Sh_\eta^\mathrm{tor})^\ad$.
 \end{itemize}
 The former is clear from the definition. 
 For the latter, note that $(\Sh_\eta^{\mathrm{tor}})^\ad\setminus \Sh_\eta^\ad=(\Sh^{\mathrm{tor}}_{\ge 1,\eta})^\ad$ is stable under generalizations \cite[Lemma 1.1.10 v)]{MR1734903}. By this fact, 
 the closure of the quasi-compact open subset $\spp^{-1}(\Sh_{v,0}^{\mathrm{tor}})$ in 
 $(\Sh_\eta^{\mathrm{tor}})^\ad$ coincides with the closure of it in $\Sh_\eta^\ad$,
 that is, $(\Sh_{\mathopen{]}\ge 1\mathclose{[}})^c$ (\cf \cite[Corollary of Theorem 1]{MR0251026}).
 Therefore $(\Sh_{\mathopen{]}\ge 1\mathclose{[}})^c$ is closed in $(\Sh_\eta^{\mathrm{tor}})^\ad$.

 By these two facts and the properness of $(\Sh^{\mathrm{tor}}_\eta)^\ad$, we have a natural homomorphism
 \begin{align*}
  &H^i_c\bigl(\spp^{-1}(\Sh_{v,0}^{\mathrm{tor}})_{\overline{\eta}},p_{m*}\mathcal{F}^\ad\bigr)\longrightarrow H^i_{(\Sh_{\mathopen{]}\ge 1\mathclose{[}})^c_{\overline{\eta}}}\bigl((\Sh^{\mathrm{tor}})^\ad_{\overline{\eta}},j_!p_{m*}\mathcal{F}^\ad\bigr)\\
  &\qquad=H^i_{(\Sh_{\mathopen{]}\ge 1\mathclose{[}})^c_{\overline{\eta}}}(\Sh^\ad_{\overline{\eta}},p_{m*}\mathcal{F}^\ad)
  =H^i_{(\Sh_{m,\mathopen{]}\ge 1\mathclose{[}})^c_{\overline{\eta}}}(\Sh^\ad_{m,\overline{\eta}},\mathcal{F}^\ad),
 \end{align*}
 where $j\colon \Sh_\eta\hooklongrightarrow \Sh^{\mathrm{tor}}_\eta$ denotes the canonical open
 immersion.
 Let us prove that this is an isomorphism. Put $\mathcal{G}=p_{m*}\mathcal{F}$.
 Consider the following distinguished triangle:
 \[
  R\Gamma_c(\Sh_{\overline{v}},R\psi \mathcal{G})\longrightarrow R\Gamma(\Sh^{\mathrm{tor}}_{\overline{v}},R\psi j_!\mathcal{G})\longrightarrow R\Gamma\bigl(\Sh^{\mathrm{tor}}_{\overline{v},\ge 1},(R\psi j_!\mathcal{G})\vert_{\Sh^{\mathrm{tor}}_{\overline{v},\ge 1}}\bigr)\yrightarrow{+1}.
 \]
 By \cite[Theorem 3.1]{MR1395723} and Lemma \ref{lem:rig-gen-fiber}, we have isomorphisms
 \begin{align*}
  R\Gamma(\Sh^{\mathrm{tor}}_{\overline{v}},R\psi j_!\mathcal{G})&\cong R\Gamma\bigl((\Sh^{\mathrm{tor}})^\ad_{\overline{\eta}},j_!\mathcal{G}^\ad\bigr),\\
  R\Gamma\bigl(\Sh^{\mathrm{tor}}_{\overline{v},\ge 1},(R\psi j_!\mathcal{G})\vert_{\Sh^{\mathrm{tor}}_{\overline{v},\ge 1}}\bigr)&\cong R\Gamma\bigl(\spp^{-1}(\Sh^{\mathrm{tor}}_{v,\ge 1})^\circ_{\overline{\eta}},j_!\mathcal{G}^\ad\bigr).
 \end{align*}
 Therefore, the natural morphism
 \[
  R\Gamma_c(\Sh_{\overline{v}},R\psi \mathcal{G})\cong R\Gamma_c\bigl(\spp^{-1}(\Sh^{\mathrm{tor}}_{v,0})_{\overline{\eta}},\mathcal{G}^\ad\bigr)\longrightarrow H^i_{(\Sh_{\mathopen{]}\ge 1\mathclose{[}})^c_{\overline{\eta}}}\bigl((\Sh^\mathrm{tor})^\ad_{\overline{\eta}},j_!\mathcal{G}^\ad\bigr)
 \]
 is an isomorphism (note that $(\Sh^\mathrm{tor}_\eta)^\ad\setminus \spp^{-1}(\Sh^{\mathrm{tor}}_{v,\ge 1})^\circ=(\Sh_{\mathopen{]}\ge 1\mathclose{[}})^c$). This completes a proof of the second isomorphism.

 By taking projective limit, we obtain natural isomorphisms
 \[
 H^i_c(\Sh_{\overline{v}},R\psi p_{m*}\mathcal{L}_\xi)\cong 
 H^i_c\bigl(\spp^{-1}(\Sh_{v,0}^{\mathrm{tor}})_{\overline{\eta}},p_{m*}\mathcal{L}_\xi^\ad\bigr)\cong
 H^i_{(\Sh_{m,\mathopen{]}\ge 1\mathclose{[}})^c_{\overline{\eta}}}(\Sh^\ad_{m,\overline{\eta}},\mathcal{L}_\xi^\ad)
 \]
 in the lemma. It is straightforward to see that the isomorphism
 \[
 \varinjlim_m H^i_c(\Sh_{\overline{v}},R\psi p_{m*}\mathcal{L_\xi})\cong 
 \varinjlim_m H^i_{(\Sh_{m,\mathopen{]}\ge 1\mathclose{[}})^c_{\overline{\eta}}}(\Sh^\ad_{m,\overline{\eta}},\mathcal{L}_\xi^\ad)
 \]
 is $G(\Q_p)$-equivariant.
\end{prf}

\begin{cor}\label{cor:exact-seq}
 We have the following long exact sequence of $G(\Q_p)$-modules:
 \begin{align*}
  \cdots&\longrightarrow H^i_c(\Sh_{\infty,K^p,\overline{v}},R\psi\mathcal{L}_\xi)
  \longrightarrow H^i(\Sh_{\infty,K^p,\overline{\eta}},\mathcal{L}_\xi)
  \longrightarrow \varinjlim_m H^i(\Sh_{m,\mathopen{]}\ge 1\mathclose{[},\overline{\eta}},\mathcal{L}^\ad_\xi)\\
  &\longrightarrow H^{i+1}_c(\Sh_{\infty,K^p,\overline{v}},R\psi\mathcal{L}_\xi)\longrightarrow \cdots,
 \end{align*}
 where the action of $G(\Q_p)$ on $\varinjlim_m H^i(\Sh_{m,\mathopen{]}\ge 1\mathclose{[},\overline{\eta}},\mathcal{L}^\ad_\xi)$ comes from Proposition \ref{prop:torus-rank-Hecke}.
\end{cor}

\begin{prf}
 Write $\mathcal{L}_\xi=(\mathcal{F}_n)_n\otimes \overline{\Q}_\ell$, where $\mathcal{F}_n$ is
 a locally constant constructible $\ell^n$-torsion sheaf.
 By the proof of Lemma \ref{lem:nearby-rigid}, we have a natural isomorphism
 \[
  H^i_c(\Sh_{K^p,\overline{v}},R\psi p_{m*}\mathcal{F}_n)\cong H^i_{(\Sh_{m,\mathopen{]}\ge 1\mathclose{[}})^c_{\overline{\eta}}}(\Sh^\ad_{m,\overline{\eta}},\mathcal{F}_n^\ad).
 \]
 On the other hand, by the comparison isomorphism \cite[Theorem 3.8.1]{MR1734903}, we have
 \[
  H^i(\Sh_{m,K^p,\overline{\eta}},\mathcal{F}_n)\cong H^i(\Sh_{m,\overline{\eta}}^\ad,\mathcal{F}^\ad_n).
 \]
 Therefore we have the following long exact sequence:
  \begin{align*}
  \cdots&\longrightarrow H^i_c(\Sh_{K^p,\overline{v}},R\psi p_{m*}\mathcal{F}_n)
  \longrightarrow H^i(\Sh_{m,K^p,\overline{\eta}},\mathcal{F}_n)
  \longrightarrow H^i(\Sh_{m,\mathopen{]}\ge 1\mathclose{[},\overline{\eta}},\mathcal{F}^\ad_n)\\
  &\longrightarrow H^{i+1}_c(\Sh_{K^p,\overline{v}},R\psi p_{m*}\mathcal{F}_n)\longrightarrow \cdots.
  \end{align*}
 On the other hand, clearly $H^i_c(\Sh_{K^p,\overline{v}},R\psi p_{m*}\mathcal{F}_n)$ and 
 $H^i(\Sh_{m,K^p,\overline{\eta}},\mathcal{F}_n)$ are finitely generated $\Z/\ell^n\Z$-modules for
 every $i$, and thus $H^i_c(\Sh_{m,\mathopen{]}\ge 1\mathclose{[},\overline{\eta}},\mathcal{F}^\ad_n)$
 is also finitely generated. Therefore the above exact sequence of projective systems satisfies
 the Mittag-Leffler condition, and the projective limit preserves exactness.
 By taking projective limit with respect to $n$, tensoring $\overline{\Q}_\ell$, and taking inductive limit
 with respect to $m$, we obtain the desired exact sequence.
\end{prf}

\begin{lem}\label{lem:coh-finiteness}
 Let $\mathcal{F}$ be a constructible $\ell^n$-torsion sheaf on $\Sh_{m,\overline{\eta}}$.
 Then, for every $r\ge 1$, $H^i(\Sh_{m,\mathopen{]}\ge r\mathclose{[},\overline{\eta}},\mathcal{F}^\ad)$
 is a finitely generated $\Z/\ell^n\Z$-module.
\end{lem}

\begin{prf}
 Let $j\colon \Sh_{\eta}\hooklongrightarrow \Sh^{\mathrm{tor}}_{\eta}$ be the canonical open immersion.
 Then, by \cite[Theorem 3.7.2, Theorem 3.8.1]{MR1734903}, \cite[Theorem 3.1]{MR1395723} and Lemma \ref{lem:rig-gen-fiber}, we have
 \begin{align*}
  H^i(\Sh_{m,\mathopen{]}\ge r\mathclose{[},\overline{\eta}},\mathcal{F}^\ad)
  &=H^i(\Sh_{\mathopen{]}\ge r\mathclose{[},\overline{\eta}},p_{m*}\mathcal{F}^\ad)
 =H^i\bigl(\spp^{-1}(\Sh^{\mathrm{tor}}_{v,\ge r})^\circ_{\overline{\eta}},Rj_*p_{m*}\mathcal{F}^\ad\bigr)\\
  &\cong H^i\bigl(\spp^{-1}(\Sh^{\mathrm{tor}}_{v,\ge r})^\circ_{\overline{\eta}},(Rj_*p_{m*}\mathcal{F})^\ad\bigr)\\
  & \cong H^i\bigl(\Sh^{\mathrm{tor}}_{\overline{v},\ge r},(R\psi Rj_*p_{m*}\mathcal{F})\vert_{\Sh^{\mathrm{tor}}_{\overline{v},\ge r}}\bigr).
 \end{align*}
 This is a finitely generated $\Z/\ell^n\Z$-module, as $(R\psi Rj_*p_{m*}\mathcal{F})\vert_{\Sh^{\mathrm{tor}}_{\overline{v},\ge r}}$ is a constructible $\Z/\ell^n\Z$-complex on $\Sh^{\mathrm{tor}}_{\overline{v},\ge r}$.
\end{prf}

\begin{cor}\label{cor:exact-sequence-2}
 Let $r\ge 1$ be an integer. Put $V^i_{\ge r}=\varinjlim_m H^i(\Sh_{m,\mathopen{]}\ge r\mathclose{[},\overline{\eta}},\mathcal{L}_\xi^\ad)$ and $V^i_{r}=\varinjlim_m H^i_{\Sh_{m,\mathopen{]}r\mathclose{[},\overline{\eta}}}(\Sh_{m,\mathopen{]}\ge r\mathclose{[},\overline{\eta}},\mathcal{L}_\xi^\ad)$.
 \begin{enumerate}
  \item The group $G(\Q_p)$ naturally acts on $V^i_{\ge r}$ and $V^i_{r}$, 
	and these are admissible $G(\Q_p)$-representations.
  \item We have the following long exact sequence of $G(\Q_p)$-modules:
	\[
	 \cdots\longrightarrow V^i_{r}\longrightarrow V^i_{\ge r}\longrightarrow V^i_{\ge r+1}
	\longrightarrow V^{i+1}_{r}\longrightarrow \cdots.
	\]
 \end{enumerate}
\end{cor}

\begin{prf}
 First let us construct the long exact sequence in ii). As in the proof of Corollary \ref{cor:exact-seq},
 write $\mathcal{L}_\xi=(\mathcal{F}_n)_n\otimes\overline{\Q}_\ell$. Obviously we have the following
 long exact sequence:
 \begin{align*}
  \cdots&\longrightarrow H^i_{\Sh_{m,\mathopen{]}r\mathclose{[},\overline{\eta}}}(\Sh_{m,\mathopen{]}\ge r\mathclose{[},\overline{\eta}},\mathcal{F}_n^\ad)\longrightarrow H^i(\Sh_{m,\mathopen{]}\ge r\mathclose{[},\overline{\eta}},\mathcal{F}_n^\ad)\longrightarrow H^i(\Sh_{m,\mathopen{]}\ge r+1\mathclose{[},\overline{\eta}},\mathcal{F}_n^\ad)\\
  &\longrightarrow H^{i+1}_{\Sh_{m,\mathopen{]}r\mathclose{[},\overline{\eta}}}(\Sh_{m,\mathopen{]}\ge r\mathclose{[},\overline{\eta}},\mathcal{F}_n^\ad)\longrightarrow\cdots.
 \end{align*}
 Lemma \ref{lem:coh-finiteness} tells us that every term of this long exact sequence is a finitely generated $\Z/\ell^n\Z$-module. Therefore, as in the proof of \ref{cor:exact-seq}, we obtain the long exact sequence
 \[
 \cdots\longrightarrow V^i_{r}\longrightarrow V^i_{\ge r}\longrightarrow V^i_{\ge r+1}
 \longrightarrow V^{i+1}_{r}\longrightarrow \cdots
 \]
 by taking projective limit, tensoring $\overline{\Q}_\ell$ and taking inductive limit.

 By Proposition \ref{prop:torus-rank-Hecke}, the group $G(\Q_p)$ acts on $V^i_{\ge r}$ and $V^i_{r}$.
 Clearly they are smooth $G(\Q_p)$-representations, and the long exact sequence above is $G(\Q_p)$-equivariant.
 Moreover, as in the proof of Lemma \ref{lem:nearby-gp-action},
 we can easily see that for an integer $m\ge 1$,
 \[
  (V^i_{\ge r})^{K_m}=H^i(\Sh_{m,\mathopen{]}\ge r\mathclose{[},\overline{\eta}},\mathcal{L}_\xi^\ad),\quad
  (V^i_{r})^{K_m}=H^i_{\Sh_{m,\mathopen{]}r\mathclose{[},\overline{\eta}}}(\Sh_{m,\mathopen{]}\ge r\mathclose{[},\overline{\eta}},\mathcal{L}_\xi^\ad).
 \]
 Since $H^i(\Sh_{m,\mathopen{]}\ge r\mathclose{[},\overline{\eta}},\mathcal{F}_n^\ad)$ and
 $H^i_{\Sh_{m,\mathopen{]}r\mathclose{[},\overline{\eta}}}(\Sh_{m,\mathopen{]}\ge r\mathclose{[},\overline{\eta}},\mathcal{F}_n^\ad)$ are finitely generated $\Z/\ell^n\Z$-modules for every $n\ge 1$, 
 $H^i(\Sh_{m,\mathopen{]}\ge r\mathclose{[},\overline{\eta}},\mathcal{L}_\xi^\ad)$ and 
 $H^i_{\Sh_{m,\mathopen{]}r\mathclose{[},\overline{\eta}}}(\Sh_{m,\mathopen{]}\ge r\mathclose{[},\overline{\eta}},\mathcal{L}_\xi^\ad)$ are finite-dimensional $\overline{\Q}_\ell$-vector spaces.
 Therefore we conclude that $V^i_{\ge r}$ and $V^i_{r}$ are admissible $G(\Q_p)$-representations.
\end{prf}

\begin{prop}\label{prop:parabolically-induced}
 Let $P_1,\ldots,P_k$ be the parabolic subgroups of $G_{\Q_p}$ as in Proposition \ref{prop:S_m-parabolic} ii).
 For an integer $j$ with $1\le j\le k$, let $I(j)\in\mathcal{S}_\infty$ be
 the image of $[\id]\in G(\Q_p)/P_j(\Q_p)\subset \coprod_{i=1}^k G(\Q_p)/P_i(\Q_p)$ under the isomorphism
 $\coprod_{i=1}^k G(\Q_p)/P_i(\Q_p)\cong \mathcal{S}_\infty$ constructed in 
 Proposition \ref{prop:S_m-parabolic} ii). For $m\ge 1$, denote by $I(j)_m$ the image of $I(j)$ under
 the natural map $\mathcal{S}_\infty\longrightarrow \mathcal{S}_m$.
 Put $W^i_{r,j}=\varinjlim_m H^i_{\Sh_{m,\mathopen{]}r\mathclose{[},I(j)_m,\overline{\eta}}}(\Sh_{m,\mathopen{]}\ge r\mathclose{[},\overline{\eta}},\mathcal{L}_\xi^\ad)$.
\begin{enumerate}
 \item For every $j$ with $1\le j\le k$, we have a natural smooth action of $P_j(\Q_p)$ on $W^i_{r,j}$.
 \item We have a natural $G(\Q_p)$-equivariant isomorphism 
	$V^i_r\cong \bigoplus_{j}\Ind_{P_j(\Q_p)}^{G(\Q_p)}W^i_{r,j}$.
 \item The $G(\Q_p)$-representation $V^i_r$ is non-cuspidal.
 \end{enumerate}
\end{prop}

\begin{prf}
 First note that $\Sh_{m+1,\mathopen{]}r\mathclose{[},I(j)_{m+1}}$ is open and closed
 in $p_{m+1,m}^{-1}(\Sh_{m,\mathopen{]}r\mathclose{[},I(j)_m})$, where
 $p_{m+1,m}\colon \Sh_{m+1,\eta}\longrightarrow \Sh_{m,\eta}$ denotes the natural transition map;
 it follows from Proposition \ref{prop:partition-I} and Proposition \ref{prop:partition-I-Hecke}.
 Therefore the transition map
 \[
  H^i_{\Sh_{m,\mathopen{]}r\mathclose{[},I(j)_m,\overline{\eta}}}(\Sh_{m,\mathopen{]}\ge r\mathclose{[},\overline{\eta}},\mathcal{L}_\xi^\ad)\longrightarrow H^i_{\Sh_{m+1,\mathopen{]}r\mathclose{[},I(j)_{m+1},\overline{\eta}}}(\Sh_{m+1,\mathopen{]}\ge r\mathclose{[},\overline{\eta}},\mathcal{L}_\xi^\ad)
 \]
 of $W_{r,j}^i$ is induced.

 Now i) is clear from Proposition \ref{prop:partition-I-Hecke}. Let us prove ii).
 We follow the proof of \cite[Proposition 5.20]{RZ-GSp4}. By the Frobenius reciprocity, 
 we have a homomorphism of $G(\Q_p)$-modules
 \[
 \bigoplus_{j=1}^k\Ind_{P_j(\Q_p)}^{G(\Q_p)}W^i_{r,j}\longrightarrow V^i_r.
 \]
 We shall observe that this is bijective. For an integer $m\ge 1$, we have
 \begin{align*}
  H^i_{\Sh_{m,\mathopen{]}r\mathclose{[},\overline{\eta}}}(\Sh_{m,\mathopen{]}\ge r\mathclose{[},\overline{\eta}},\mathcal{L}_\xi^\ad)&\stackrel{(1)}{=}\bigoplus_{I\in\mathcal{S}_m} H^i_{\Sh_{m,\mathopen{]}r\mathclose{[},I,\overline{\eta}}}(\Sh_{m,\mathopen{]}\ge r\mathclose{[},\overline{\eta}},\mathcal{L}_\xi^\ad)\\
  &\stackrel{(2)}{=}\bigoplus_{j=1}^k\bigoplus_{g\in K_m\backslash G(\Q_p)/P_j(\Q_p)} H^i_{\Sh_{m,\mathopen{]}r\mathclose{[},g^{-1}I(j)_m,\overline{\eta}}}(\Sh_{m,\mathopen{]}\ge r\mathclose{[},\overline{\eta}},\mathcal{L}_\xi^\ad)\\
  &\stackrel{(3)}{=}\bigoplus_{j=1}^k\bigoplus_{g\in K_m\backslash K_0/P_j(\Q_p)\cap K_0} H^i_{\Sh_{m,\mathopen{]}r\mathclose{[},g^{-1}I(j)_m,\overline{\eta}}}(\Sh_{m,\mathopen{]}\ge r\mathclose{[},\overline{\eta}},\mathcal{L}_\xi^\ad)\\
  &\stackrel{(4)}{\cong} \bigoplus_{j=1}^k\Ind_{(P_j(\Q_p)\cap K_0)/(P_j(\Q_p)\cap K_m)}^{K_0/K_m}H^i_{\Sh_{m,\mathopen{]}r\mathclose{[},I(j)_m,\overline{\eta}}}(\Sh_{m,\mathopen{]}\ge r\mathclose{[},\overline{\eta}},\mathcal{L}_\xi^\ad).
 \end{align*}
 Here (1) follows from Proposition \ref{prop:partition-I}, (2) from Proposition \ref{prop:S_m-parabolic},
 (3) from the Iwasawa decomposition $G(\Q_p)=P_j(\Q_p)K_0$.
 The isomorphism (4) is a consequence of Proposition \ref{prop:partition-I-Hecke} and \cite[Lemme 13.2]{MR1719811}.
 By taking the inductive limit, we obtain $K_0$-isomorphisms
 \[
  V^i_{r}\cong \bigoplus_{j=1}^k \Ind_{P_j(\Q_p)\cap K_0}^{K_0}W^i_{r,j}\yleftarrow{\cong}\Ind_{P_j(\Q_p)}^{G(\Q_p)}W^i_{r,j}
 \]
 (the second map is an isomorphism by the Iwasawa decomposition).
 By the proof of \cite[Lemme 13.2]{MR1719811},
 it is easy to see that the first isomorphism above is nothing but the $K_0$-homomorphism obtained by
 the Frobenius reciprocity for $P_j(\Q_p)\cap K_0\subset K_0$. Therefore the composite of the two isomorphisms
 above coincides with the $G(\Q_p)$-homomorphism introduced at the beginning of our proof of ii).
 Thus we conclude the proof of ii).

 Finally consider iii). By ii), we have only to prove that the unipotent radical of $P_j(\Q_p)$
 acts trivially on $W^i_{r,j}$. By \cite[Lemme 13.2.3]{MR1719811}, it suffices to prove that
 $W^i_{r,j}$ is an admissible $P_j(\Q_p)$-representation.
 For an integer $m\ge 1$, $(W^i_{r,j})^{P_j(\Q_p)\cap K_m}$ is a subspace of
 $(\Ind_{P_j(\Q_p)}^{G(\Q_p)}W^i_{r,j})^{K_m}$.
 By ii) and Corollary \ref{cor:exact-sequence-2} i),
 $(\Ind_{P_j(\Q_p)}^{G(\Q_p)}W^i_{r,j})^{K_m}\subset (\bigoplus_{j=1}^k\Ind_{P_j(\Q_p)}^{G(\Q_p)}W^i_{r,j})^{K_m}\cong (V^i_r)^{K_m}$ is a finite-dimensional $\overline{\Q}_\ell$-vector space.
 Hence $W^i_{r,j}$ is an admissible $P_j(\Q_p)$-representation, as desired.
\end{prf}

\begin{prf}[of Theorem \ref{thm:main-thm}]
 By Remark \ref{rem:cpt-ordinary} and Corollary \ref{cor:exact-seq}, it suffices to show that
 $V^i_{\ge 1}=\varinjlim_m H^i(\Sh_{m,\mathopen{]}\ge 1\mathclose{[},\overline{\eta}},\mathcal{L}^\ad_\xi)$ is
 non-cuspidal for every $i$. By Corollary \ref{cor:exact-sequence-2} and Proposition \ref{prop:parabolically-induced}, the non-cuspidality of $V^i_{\ge r}$ is equivalent to that of $V^i_{\ge r+1}$. If $r$ is large enough, $V^i_{\ge r}=0$. Hence we can conclude that $V^i_{\ge r}$ is non-cuspidal for every $r\ge 1$, by the descending induction on $r$.
\end{prf}

\subsection{Torsion coefficients}
Since our method of proving Theorem \ref{thm:main-thm} is totally geometric,
we may also obtain an analogous result for $\ell$-torsion coefficients.
For simplicity, we will only consider a constant coefficient $\overline{\F}_\ell$.
As in Section \ref{subsec:cohomology-notation}, put
\begin{align*}
 H^i_c(\Sh_{\infty,K^p,\overline{\eta}},\overline{\F}_\ell)&=\varinjlim_{m}H^i_c(\Sh_{m,K^p,\eta}\otimes_{E_v}\overline{E}_v,\overline{\F}_\ell),\\
 H^i_c(\Sh_{\infty,K^p,\overline{v}},R\psi\overline{\F}_\ell)&=\varinjlim_{m}H^i_c\bigl(\Sh_{K^p,\overline{v}},R\psi (p_{m*}\overline{\F}_\ell)\bigr).
\end{align*}
They are naturally endowed with actions of $G(\Q_p)\times \Gal(\overline{E}_v/E_v)$.
They are admissible/continuous $G(\Q_p)\times \Gal(\overline{E}_v/E_v)$-representations; note that if $m\ge 1$ we have
\begin{align*}
 H^i_c(\Sh_{\infty,K^p,\overline{\eta}},\overline{\F}_\ell)^{K_m}&=H^i_c(\Sh_{m,K^p,\overline{\eta}},\overline{\F}_\ell),\\
 H^i_c(\Sh_{\infty,K^p,\overline{v}},R\psi\overline{\F}_\ell)^{K_m}&=H^i_c\bigl(\Sh_{K^p,\overline{v}},R\psi (p_{m*}\overline{\F}_\ell)\bigr),
\end{align*}
since $K_m$ is a pro-$p$ group (\cf \cite[Proposition 2.5]{non-cusp}).

The following theorem can be proven in exactly the same way as Theorem \ref{thm:main-thm}.

\begin{thm}\label{thm:main-thm-mod-l}
  The kernel and the cokernel of the canonical homomorphism
 \[
  H^i_c(\Sh_{\infty,K^p,\overline{v}},R\psi\overline{\F}_\ell)\longrightarrow H^i_c(\Sh_{\infty,K^p,\overline{\eta}},\overline{\F}_\ell)
 \]
 have no supercuspidal subquotient of $G(\Q_p)$.
 For the definition of supercuspidal representations over $\overline{\F}_\ell$, see \cite[II.2.5]{MR1395151}.
\end{thm}

\section{Applications}
In this section, we give a very simple application of our main result.
Proofs in this section are rather sketchy, since the technique is more or less well-known.
More detailed and general study will be given elsewhere. 

Here we consider the Shimura variety for $\mathit{GU}(1,n-1)$ over $\Q$.
Let $F$ be an imaginary quadratic extention of $\Q$ and $\Spl_{F/\Q}$ the set of rational primes over which $F/\Q$ splits. Fix an integer $n\ge 2$ and a prime $p\in \Spl_{F/\Q}$.
Consider the integral PEL datum
$(B,\mathcal{O}_B,*,V,L,\langle\ ,\ \rangle,h)$ as follows:
\begin{itemize}
 \item $B=F$, $\mathcal{O}_B=\mathcal{O}_F$ and $*$ is the unique non-trivial element of $\Gal(F/\Q)$.
 \item $V=F^n$ and $L=\mathcal{O}_F^n$.
 \item $\langle\ ,\ \rangle\colon V\times V\longrightarrow \Q$ be an alternating pairing
       satisfying the following conditions:
       \begin{itemize}
	\item[$\bullet$] $\langle x,y\rangle\in \Z$ for every $x,y\in L$, 
	\item[$\bullet$] $\langle bx,y\rangle=\langle x,b^*y\rangle$ for every $x,y\in V$ and $b\in F$, and
	\item[$\bullet$] $G_{\R}\cong \mathit{GU}(1,n-1)$ (for the definition of $G$, see Section \ref{subsec:notation-Sh-var}).
       \end{itemize}
 \item $h\colon \C\longrightarrow \End_F(V)\otimes_\Q\R\cong M_n(\C)$ is given by
       $z\longmapsto \begin{pmatrix}z&0\\0&\overline{z}I_{n-1}\end{pmatrix}$.
\end{itemize}
Put $\Sigma=\{p\in\Spl_{F/\Q}\mid L_p=L_p^\perp\}$. Our integral Shimura datum is unramified at every $p\in\Sigma$.
Moreover, for such $p$, $G_{\Q_p}$ is isomorphic to $\GL_n(\Q_p)\times \GL_1(\Q_p)$
(\cf \cite[\S 1.2.3]{MR2074714}).
In this case, the reflex field $E$ is equal to $F$. To a neat compact open subgroup $K$ of $G(\widehat{\Z})$,
we can attach the Shimura variety $\Sh_K$, which is not proper over $\Spec F$. 
If $K=K_0K^p$ for some compact open subgroup
$K^p$ of $G(\widehat{\Z}^p)$, we have $\Sh_K=\Sh_{K^p}\otimes_{\mathcal{O}_{F,(p)}}F$, where
$\Sh_{K^p}$ is the moduli space introduced in Section \ref{subsec:notation-Sh-var}.

Let us fix a prime number $\ell$. Put 
\[
 H_c^i(\Sh,\overline{\Q}_\ell)=\varinjlim_{K}H_c^i(\Sh_K\otimes_F\overline{F},\overline{\Q}_\ell).
\]
It is an admissible/continuous $G(\A^\infty)\times\Gal(\overline{F}/F)$-representation over $\overline{\Q}_\ell$.

\begin{thm}\label{thm:non-cusp}
 Let $\Pi$ be an irreducible admissible representation of $G(\A^\infty)$ over $\overline{\Q}_\ell$.
 Assume that there exists a prime $p\in \Sigma$ such that $\Pi_p$ is a supercuspidal 
 representation of $G(\Q_p)$.
 Then $H_c^i(\Sh,\overline{\Q}_\ell)[\Pi]=0$ unless $i=n-1$.
\end{thm}

\begin{rem}
 For proper Shimura varieties, an analogous result is known (\cite{MR1114211}, \cite[Corollary IV.2.7]{MR1876802}).
 It would be possible to give an ``automorphic'' proof of Theorem \ref{thm:non-cusp}
 by using results in \cite{MR2567740}. However, the authors think that our proof,
 consisting of purely local arguments, is simpler and has some importance.
\end{rem}

\begin{prf}
 Let $\ell'$ be another prime number and fix an isomorphism of fields $\iota\colon\overline{\Q}_\ell\cong \overline{\Q}_{\ell'}$. Then $\iota$ induces an isomorphism $H_c^i(\Sh,\overline{\Q}_\ell)[\Pi]\cong H_c^i(\Sh,\overline{\Q}_{\ell'})[\iota \Pi]$, where $\iota \Pi$ is the representation of $G(\A^\infty)$ over $\overline{\Q}_{\ell'}$
 induced by $\Pi$ and $\iota$. It is easy to observe that $\Pi_p$ is supercuspidal if and only if
 $(\iota\Pi)_p$ is supercuspidal. Therefore, we can change our $\ell$ freely, and thus we can assume that
 there exists a prime $p\in \Sigma\setminus \{\ell\}$ such that $\Pi_p$ is supercuspidal.  
 Fix such $p$ and take a place $v$ of $F$ lying over $p$. Then, for an integer $m\ge 0$ and a neat compact open
 subgroup $K^p$ of $G(\widehat{\Z}^p)$, $\Sh_{K_mK^p}\otimes_FF_v$ is isomorphic to $\Sh_{m,K^p,\eta}$
 introduced in Section \ref{subsec:notation-Sh-var}.
 Therefore we have an isomorphism
 \[
  H_c^i(\Sh,\overline{\Q}_\ell)\cong \varinjlim_{m,K^p}H_c^i(\Sh_{m,K^p,\overline{\eta}},\overline{\Q}_\ell)
 =\varinjlim_{K^p}H_c^i(\Sh_{\infty,K^p,\overline{\eta}},\overline{\Q}_\ell).
 \]
 Thus it suffices to show that $H_c^i(\Sh_{\infty,K^p,\overline{\eta}},\overline{\Q}_\ell)[\pi]=0$
 for a supercuspidal representation $\pi$ of $G(\Q_p)$, a neat compact open subgroup $K^p$, 
 and an integer $i\neq n-1$.
 By Theorem \ref{thm:main-thm}, it is equivalent to showing that
 $H_c^i(\Sh_{\infty,K^p,\overline{v}},R\psi\overline{\Q}_\ell)[\pi]=0$.

 For an integer $h\ge 0$, let $\Sh^{[h]}_{K^p,v}$ be the reduced closed subscheme of $\Sh_{K^p,v}$
 consisting of points $x$ such that the \'etale rank of $\mathcal{A}_x[v^\infty]$ is less than
 or equal to $h$ (\cf \cite[p.~111]{MR1876802}); 
 recall that $\mathcal{A}$ denotes the universal abelian scheme over $\Sh_{K^p}$.
 Put $\Sh^{(h)}_{K^p,v}=\Sh^{[h]}_{K^p,v}\setminus \Sh^{[h-1]}_{K^p,v}$.
 Our proof of the theorem is divided into the subsequent two lemmas.
\end{prf} 

\begin{lem}
 For every supercuspidal representation $\pi$ of $G(\Q_p)$, 
 we have 
 \[
  H_c^i(\Sh_{\infty,K^p,\overline{v}},R\psi\overline{\Q}_\ell)[\pi]
 =\Bigl(\varinjlim_{m}H_c^i\bigl(\Sh^{[0]}_{K^p,\overline{v}},(R\psi p_{m*}\overline{\Q}_\ell)\vert_{\Sh^{[0]}_{K^p,\overline{v}}}\bigr)\Bigr)[\pi].
 \]
\end{lem}

\begin{prf}
 First recall that our Shimura variety $\Sh_{m,K^p,\eta}$ has a good integral model over $\mathcal{O}_v$. 
 For an integer $m\ge 0$, consider the functor from the category of $\mathcal{O}_v$-schemes to the category of
 sets, that associates $S$ to the set of isomorphism classes of 6-tuples $(A,i,\lambda,\eta^p,\eta_v,\eta_{p,0})$,
 where
 \begin{itemize}
  \item $[(A,i,\lambda,\eta^p)]\in \Sh_{K^p}(S)$,
  \item $\eta_v\colon L\otimes_\Z(v^{-m}\mathcal{O}_v/\mathcal{O}_v)\longrightarrow A[v^m]$ be a Drinfeld $v^m$-level structure
	(\cf \cite[II.2]{MR1876802}), and
  \item $\eta_{p,0}\colon p^{-m}\Z/\Z\longrightarrow \mu_{p^m,S}$ be a Drinfeld $p^m$-level structure.
 \end{itemize}
 Then it is easy to see that this functor is represented by a scheme $\Sh_{m,K^p}$
 which is finite over $\Sh_{K^p}$. Moreover the generic fiber of $\Sh_{m,K^p}$ can be naturally identified with 
 $\Sh_{m,K^p,\eta}$ (\cf the moduli problem $\mathfrak{X}'_U$ introduced in
 \cite[p.~92]{MR1876802}). As in \cite[III.4]{MR1876802}, we can extend the Hecke action of $G(\Q_p)$ on
 $(\Sh_{m,K^p,\eta})_{m\ge 0}$ to the tower $(\Sh_{m,K^p})_{m\ge 0}$. We have a $G(\Q_p)$-equivariant isomorphism
 \[
  H_c^i(\Sh_{\infty,K^p,\overline{v}},R\psi\overline{\Q}_\ell)\cong \varinjlim_{m}H_c^i\bigl(\Sh_{m,K^p,\overline{v}},R\psi \overline{\Q}_\ell\bigr).
 \]
 Let us denote by $\Sh^{[h]}_{m,K^p,v}$ (resp.\ $\Sh^{(h)}_{m,K^p,v}$) the inverse image of
 $\Sh^{[h]}_{K^p,v}$ (resp.\ $\Sh^{(h)}_{K^p,v}$) under $\Sh_{m,K^p}\longrightarrow \Sh_{K^p}$.
 For an integer $h\ge 0$, it is easy to observe that
 \begin{align*}
  \varinjlim_{m}H_c^i\bigl(\Sh^{[h]}_{K^p,\overline{v}},(R\psi p_{m*}\overline{\Q}_\ell)\vert_{\Sh^{[h]}_{K^p,\overline{v}}}\bigr)&\cong \varinjlim_{m}H_c^i\bigl(\Sh^{[h]}_{m,K^p,\overline{v}},(R\psi \overline{\Q}_\ell)\vert_{\Sh^{[h]}_{m,K^p,\overline{v}}}\bigr),\\
 \varinjlim_{m}H_c^i\bigl(\Sh^{(h)}_{K^p,\overline{v}},(R\psi p_{m*}\overline{\Q}_\ell)\vert_{\Sh^{(h)}_{K^p,\overline{v}}}\bigr)&\cong \varinjlim_{m}H_c^i\bigl(\Sh^{(h)}_{m,K^p,\overline{v}},(R\psi \overline{\Q}_\ell)\vert_{\Sh^{(h)}_{m,K^p,\overline{v}}}\bigr)
 \end{align*}
 and that they are admissible $G(\Q_p)$-representations. Moreover, by considering the kernel of 
 the universal Drinfeld $v^m$-level structure $\eta_v^{\mathrm{univ}}$, we can decompose $\Sh_{m,K^p,v}^{(h)}$ into finitely many
 open and closed subsets indexed by the set consisting of direct summands of $L\otimes_\Z(v^{-m}\mathcal{O}_v/\mathcal{O}_v)$ with
 rank $n-h$ (\cf \cite[D\'efinition 10.4.1, Proposition 10.4.2]{MR1719811} and \cite[Definition 5.1, Lemma 5.3]{RZ-GSp4}).
 Using this partition, we can prove that the $G(\Q_p)$-representation 
 \[
  \varinjlim_{m}H_c^i\bigl(\Sh^{(h)}_{m,K^p,\overline{v}},(R\psi \overline{\Q}_\ell)\vert_{\Sh^{(h)}_{m,K^p,\overline{v}}}\bigr)
 \]
 is parabolically induced from a proper parabolic subgroup of $G(\Q_p)$.
 Therefore, by the same argument as in the proof of Theorem \ref{thm:main-thm},
 we can conclude that the kernel and the cokernel of
 \[
  \varinjlim_{m}H_c^i(\Sh_{m,K^p,\overline{v}},R\psi \overline{\Q}_\ell)\longrightarrow \varinjlim_{m}H_c^i\bigl(\Sh^{[0]}_{m,K^p,\overline{v}},(R\psi \overline{\Q}_\ell)\vert_{\Sh^{[0]}_{m,K^p,\overline{v}}}\bigr)
 \]
 are non-cuspidal. This completes the proof of the lemma.
\end{prf}

\begin{lem}\label{lem:basic-noncusp}
 Let $\pi$ be a supercuspidal representation of $G(\Q_p)$.
 If $i\neq n-1$, we have
 \[
  \Bigl(\varinjlim_{m}H_c^i\bigl(\Sh^{[0]}_{K^p,\overline{v}},(R\psi p_{m*}\overline{\Q}_\ell)\vert_{\Sh^{[0]}_{K^p,\overline{v}}}\bigr)\Bigr)[\pi]=0.
 \]
\end{lem}

\begin{prf}
 Let $\mu_h\colon \mathbb{G}_{m,\C}\longrightarrow G_{\C}$ be the homomorphism of algebraic groups over $\C$
 defined as the composite of
 \[
  \mathbb{G}_{m,\C}\yrightarrow{z\longmapsto (z,1)}\mathbb{G}_{m,\C}\times \mathbb{G}_{m,\C}\stackrel{(*)}{\cong} (\Res_{\C/\R}\mathbb{G}_{m,\C})\otimes_\R\C\yrightarrow{h_\C}G_{\C},
 \]
 where $(*)$ is given by $\C\otimes_\R\C\yrightarrow{\cong}\C\times\C$; $a\otimes b\longmapsto (ab,a\overline{b})$.
 Fix an isomorphism of fields $\C\cong \overline{\Q}_p$ and denote by $\mu\colon \mathbb{G}_{m,\overline{\Q}_p}\longrightarrow G_{\overline{\Q}_p}$ the induced cocharacter of $G_{\overline{\Q}_p}$.
 Let $b$ be a unique basic element of $B(G_{\Q_p},\mu)$
 (for the definition of $B(G,\mu)$, we refer to \cite[\S 2.1.1]{MR2074714}),
 and denote by $\mathcal{M}$ the Rapoport-Zink space associated to the local unramified PEL datum
 $(F\otimes_{\Q}\Q_p,\mathcal{O}_F\otimes_{\Z}\Z_p,*,V_p,L_p,\langle\ ,\ \rangle,b,\mu)$
 (\cf \cite[\S 2.3.5]{MR2074714}). The Rapoport-Zink space $\mathcal{M}$ is equipped with an action of
 the group $J(\Q_p)$,
 where $J$ denotes the algebraic group over $\Q_p$ associated to $b$ (\cf \cite[Proposition 1.12]{MR1393439}).
 By \cite[\S 2.3.7.1]{MR2074714}, 
 $\mathcal{M}$ is isomorphic to $\mathcal{M}_{\mathrm{LT}}\times\Q_p^\times/\Z_p^\times$, where 
 $\mathcal{M}_{\mathrm{LT}}$ is the Lubin-Tate space for $\GL_n(\Q_p)$.
 Furthermore, $J(\Q_p)$ is isomorphic to $D^\times\times \Q_p^\times$, where $D$ denotes the central
 division algebra over $\Q_p$ with invariant $1/n$. The action of $J(\Q_p)$ on $\mathcal{M}$ is identified with
 the well-known action of $D^\times\times \Q_p^\times$ on $\mathcal{M}_{\mathrm{LT}}\times\Q_p^\times/\Z_p^\times$.

 By the $p$-adic uniformization theorem of Rapoport-Zink (\cite[Theorem 6.30]{MR1393439}, \cite[Corollaire 3.1.9]{MR2074714}), we have an isomorphism
 \[
 \coprod_{\ker^1(\Q,G)}I(\Q)\backslash \mathcal{M}\times G(\A^{\infty,p})/K^p\cong \Sh_{K^p}^{\wedge},
 \]
 where $I$ is an algebraic group over $\Q$ satisfying $I(\A^\infty)\cong J(\Q_p)\times G(\A^{\infty,p})$ and
 $\Sh_{K^p}^{\wedge}$ denotes the formal completion of $\Sh_{K^p}\otimes_{\mathcal{O}_v}W(\overline{\F}_p)$
 along $\Sh_{K^p,\overline{v}}^{[0]}$, the basic locus of $\Sh_{K^p,\overline{v}}$.
 By this isomorphism, we know that $\Sh^{[0]}_{K^p,\overline{v}}$, which coincides with $\Sh_{K^p}^{\wedge}$
 as topological spaces, consists of finitely many closed points; indeed, 
 the left hand side of the isomorphism above is a finite disjoint sum of formal schemes of the form
 $\Gamma\backslash \mathcal{M}$, where $\Gamma\subset J(\Q_p)$ is a discrete cocompact subgroup
 (\cf \cite[Lemme 3.1.7]{MR2074714}). Therefore,  by \cite[Theorem 3.1]{MR1395723},
 we have an isomorphism
 \begin{align*}
 H_c^i\bigl(\Sh^{[0]}_{K^p,\overline{v}},(R\psi p_{m*}\overline{\Q}_\ell)\vert_{\Sh^{[0]}_{K^p,\overline{v}}}\bigr)
  &=H^i\bigl(\Sh^{[0]}_{K^p,\overline{v}},(R\psi p_{m*}\overline{\Q}_\ell)\vert_{\Sh^{[0]}_{K^p,\overline{v}}}\bigr)\\
  &\cong H^i\bigl(\Sh_{m,K^p,\overline{\eta}}(b),\Q_\ell\bigr),
 \end{align*}
 where $\Sh_{m,K^p}(b)=p_m^{-1}(\spp^{-1}(\Sh_{K^p,v}^{[0]})^\circ)$.

 Now use the Hochschild-Serre spectral sequence \cite[Th\'eor\`eme 4.5.12]{MR2074714}
\begin{align*}
 &E_2^{r,s}=\varinjlim_{m}\Ext^r_{J(\Q_p)\text{-smooth}}\bigl(H^{2(n-1)-s}_c(\mathcal{M}_{K_m},\overline{\Q}_\ell)(n-1),\mathcal{A}(I)_{\mathbf{1}}^{K^p}\bigr)\\
 &\qquad\qquad\qquad\Longrightarrow \varinjlim_{m}H^{r+s}(\Sh_{m,K^p,\overline{\eta}}(b),\overline{\Q}_\ell),
\end{align*}
 where $\mathcal{A}(I)_{\mathbf{1}}$ is the space of automorphic forms on $I(\A^\infty)$
 (see \cite[D\'efinition 4.5.8]{MR2074714} for detail).
 Since $J(\Q_p)=D^\times\times \Q_p^\times$ is anisotropic modulo center, 
 it is easy to see that $E_2^{r,s}=0$ unless $r=0$. If $r=0$, we have
 \begin{align*}
  E_2^{0,s}&=\varinjlim_m\Hom_{J(\Q_p)}\bigl(H^{2(n-1)-s}_c(\mathcal{M}_\infty,\overline{\Q}_\ell)(n-1),\mathcal{A}(I)_{\mathbf{1}}^{K^p}\bigr)^{K_m},
 \end{align*}
 where we put $H^i_c(\mathcal{M}_\infty,\overline{\Q}_\ell)=\varinjlim_{m}H^i_c(\mathcal{M}_{K_m},\overline{\Q}_\ell)$.

 By \cite{non-cusp},
 the $G(\Q_p)$-representation $H^{2(n-1)-s}_c(\mathcal{M}_\infty,\overline{\Q}_\ell)(n-1)$ has
 non-zero supercuspidal part only if $s=n-1$. Indeed, for an irreducible supercuspidal representation
 $\pi=\pi_1\otimes \chi$ of $G(\Q_p)$, where $\pi_1$ is an irreducible supercuspidal representation of $\GL_n(\Q_p)$
 and $\chi$ is a character of $\GL_1(\Q_p)$, we have
 \[
 H^i_c(\mathcal{M}_\infty,\overline{\Q}_\ell)[\pi]
 =H^i_c(\mathcal{M}_{\mathrm{LT},\infty},\overline{\Q}_\ell)[\pi_1]\otimes\chi,
 \]
 as we see in \cite[p.~168]{MR2074714}. Therefore $E_2^{0,s}$ has a supercuspidal subquotient only if
 $s=n-1$.

 Hence we can conclude that
 \[
 \varinjlim_m H_c^i\bigl(\Sh^{[0]}_{K^p,\overline{v}},(R\psi p_{m*}\overline{\Q}_\ell)\vert_{\Sh^{[0]}_{K^p,\overline{v}}}\bigr)\cong \varinjlim_{m}H^i(\Sh_{m,K^p,\overline{\eta}}(b),\overline{\Q}_\ell)
 \]
 has non-zero supercuspidal part only if $i=n-1$.
\end{prf}

We also have a similar result for the torsion coefficient case.
For a neat compact open subgroup $K^p$ of $G(\widehat{\Z}^p)$, put
\[
 H_c^i(\Sh_{K^p},\overline{\F}_\ell)=\varinjlim_{m}H_c^i(\Sh_{K_mK^p}\otimes_F\overline{F},\overline{\F}_\ell).
\]
It is an admissible/continuous $G(\Q_p)\times\Gal(\overline{F}/F)$-representation over $\overline{\F}_\ell$.

\begin{thm}\label{thm:non-cusp-torsion}
 Let $p$ be a prime in $\Sigma\setminus \{\ell\}$ and $\pi$ an irreducible supercuspidal $\overline{\F}_\ell$-representation of $G(\Q_p)$.
 Then, for every neat compact open subgroup $K^p$ of $G(\widehat{\Z}^p)$, we have
 $H_c^i(\Sh_{K^p},\overline{\F}_\ell)[\pi]=0$ unless $i=n-1$.
\end{thm}

\begin{rem}
 \begin{enumerate}
  \item Theorem \ref{thm:non-cusp-torsion} for proper Shimura varieties is due to S.~W.~Shin.  
	His method, using Mantovan's formula, is slightly different from ours.
  \item Using the result in \cite{Dat-ltmodl}, it is possible to describe
	the action of $W_{\Q_p}$ on $H_c^{n-1}(\Sh_{K^p},\overline{\F}_\ell)[\pi]$
	by means of the mod-$\ell$ local Langlands correspondence.
	Such study has also been carried out by S.~W.~Shin when the Shimura variety is proper.
 \end{enumerate}
\end{rem}

\begin{prf}
 Almost all arguments in the proof of Theorem \ref{thm:non-cusp} work well.
 The only one point which should be modified is about the vanishing of the supercuspidal part of $E_2^{r,s}$
 for $(r,s)\neq (0,n-1)$ in the proof of Lemma \ref{lem:basic-noncusp}; note that an irreducible $\overline{\F}_\ell$-representation of $J(\Q_p)$,
 being supercuspidal, is not necessarily injective in the category of smooth 
 $\overline{\F}_\ell$-representations of $J(\Q_p)$ with the fixed central character.
 For this point, we can use the same argument as that by Shin, in which he uses
 the vanishing of the supercuspidal part 
 $H^i_c(\mathcal{M}_{\mathrm{LT},\infty},\overline{\F}_\ell)_{\mathrm{sc}}$ for $i\neq n-1$ 
 (\cf \cite[proof of Proposition 3.1.1, Remarque 3.1.5]{Dat-ltmodl}) and 
 the projectivity of the $D^\times$-representation 
 $H^{n-1}_c(\mathcal{M}_{\mathrm{LT},\infty},\overline{\F}_\ell)_{\mathrm{sc}}$
 (\cf \cite[\S 3.2.2, Remarque iii)]{Dat-ltmodl}).
\end{prf}

\def\cprime{$'$} \def\cprime{$'$}
\providecommand{\bysame}{\leavevmode\hbox to3em{\hrulefill}\thinspace}
\providecommand{\MR}{\relax\ifhmode\unskip\space\fi MR }
\providecommand{\MRhref}[2]{%
  \href{http://www.ams.org/mathscinet-getitem?mr=#1}{#2}
}
\providecommand{\href}[2]{#2}

\bigbreak\bigbreak

\noindent Naoki Imai\par
\noindent Research Institute for Mathematical Sciences, Kyoto University, Kyoto, 606--8502, Japan\par
\noindent E-mail address: \texttt{naoki@kurims.kyoto-u.ac.jp}

\bigbreak

\noindent Yoichi Mieda\par
\noindent Faculty of Mathematics, Kyushu University, 744 Motooka, Nishi-ku, Fukuoka, 819--0395, Japan\par
\noindent E-mail address: \texttt{mieda@math.kyushu-u.ac.jp}

\end{document}